\documentclass[12pt]{article}
\usepackage{amsmath,enumerate,amsfonts,amssymb,color,graphicx}

\def\RR{{\mathbb R}}
\def\HH{{\mathbb H}}

\def\SSphere{{\mathbb S}}

\def\Ric{{\rm Ric}} 
 
\def\diag{{\rm diag}} 
\def\muGp{{\mu_\Gamma^+}} 

\newtheorem{theorem}{Theorem}[section]
\newtheorem{proposition}[theorem]{Proposition}

\newtheorem{lemma}[theorem]{Lemma}

\newtheorem{remark}[theorem]{Remark}
\newtheorem{problem}[theorem]{Problem}

\def\bproof{\noindent{\bf Proof.\;}}
\def\eproof{\hfill$\square$\medskip}
\def\dist{{\rm dist}\,}

\def\deg{{\rm deg}\,}

\DeclareFontFamily{OT1}{rsfs}{} 
\DeclareFontShape{OT1}{rsfs}{m}{n}{ <-7> rsfs5 <7-10> rsfs7 <10-> rsfs10}{} 
\DeclareMathAlphabet{\mycal}{OT1}{rsfs}{m}{n}

\newcounter{marnote}

\begin{document}
\title{A compactness theorem for a fully nonlinear Yamabe problem under a lower Ricci curvature bound}
\author{YanYan Li \footnote{Department of Mathematics, Rutgers University.
Partially supported by NSF grant DMS-1203961.
}~ and Luc Nguyen \footnote{Department of Mathematics, Princeton University}} 

\date{December 3, 2012}
\maketitle 
\begin{abstract}
We prove compactness of  
 solutions of a fully nonlinear Yamabe problem satisfying a lower Ricci curvature bound, when the manifold is not conformally diffeomorphic to the standard sphere.  This allows us to prove the existence of solutions when the associated cone $\Gamma$ satisfies $\mu^+_\Gamma\le 1$, which includes the $\sigma_k-$Yamabe problem for $k$ not smaller than half of the dimension of the manifold.
\end{abstract}
\section{Introduction}

Let $(M,g)$ be a compact smooth Riemannian manifold of dimension $n\ge 3$. Throughout the paper $M$ will always be connected. Let $A_g$ be
the Schouten tensor of $g$:
\[
A_g = \frac{1}{n-2}\Big(\Ric_g- \frac{1}{2(n-1)}\,R_g\,
g\Big), 
\] 
where $\Ric_g$ and $R_g$ are respectively the Ricci curvature and the scalar curvature of $g$. Let $\lambda(A_g)=(\lambda_1, \cdots, \lambda_n)$
denote the eigenvalues of
$A_g$ with respect to $g$, and 
\begin{align}  
&\Gamma\subset \Bbb R^n \mbox{ be an open convex symmetric
cone with vertex at  
the origin,} 
	\label{2} \\
&\{\lambda \in \Bbb R^n | \lambda_i > 0, 1 \leq i \leq n\} \subset
\Gamma \subset \{\lambda \in \Bbb R^n | \lambda_1 + \ldots + \lambda_n >
0\},
	\label{3}\\
&f\in C^\infty (\Gamma) \cap C^0 (\overline{\Gamma}) \mbox{ be concave,
homogeneous of degree one,
symmetric in } \lambda_i,
	\label{5}\\
&f>0\ \mbox{in}\ \Gamma,
\quad f = 0 \mbox{ on } \partial\Gamma ; \quad
f_{\lambda_i} > 0 \ \mbox{in
} \Gamma   \  \forall 1 \leq i \leq n.
\label{6}
\end{align}
For a positive function $u$, let $g_u$ denote
 the metric $u^{\frac{4}{n-2}}g$. Note that the Schouten tensor of $g_u$ is given by 
\[ 
A_{g_u} = -\frac{2}{n-2} u^{-1}\,\nabla_g^2 u + \frac{2n}{(n-2)^2}\,u^{-2} du \otimes du - \frac{2}{(n-2)^2}u^{-2}\,|du|_g^2\,g + A_g. 
\] 

\begin{problem}\label{MQuestClosed} 
Let $(f, \Gamma)$ satisfy
(\ref{2})-(\ref{6}),
  and $(M,g)$ be a compact, smooth
Riemannian manifold of dimension $n\ge 3$ satisfying
$\lambda(A_g)\in \Gamma$ on $M$.  Is there a
smooth positive function $u$ on $M$ such that
\begin{equation}
f(\lambda(A_{g_u})) =1, \quad \lambda (A_{g_u}) \in \Gamma,
\quad \mbox{ on } M?
\label{9}
\end{equation}
\end{problem} 
 
A closely related problem is on the compactness of the solution set.

\begin{problem}\label{MQuest-extra}
Under the hypotheses of Problem \ref{MQuestClosed}, and
assuming that $(M,g)$ is not conformally equivalent to the
standard sphere,  do all smooth positive 
solutions  of (\ref{9}) satisfy
\[
\|u\|_{ C^3(M) }+\|u^{-1}\|_{ C^3(M) }\le C,
\]
for some constant $C$ depending only on $(M,g)$
and $(f,\Gamma)$?
\end{problem}
 
Equation (\ref{9}) is a second
order fully nonlinear elliptic equation of $u$.
For $1\leq k \leq n,$ let $\sigma_k(\lambda) = \sum_{1\leq i_1 <\cdots <
i_k\leq n}\lambda_{i_1}\cdots \lambda_{i_k}, \lambda = (\lambda_1,
\cdots, \lambda_n) \in \Bbb R^n$, denote the $k$-th elementary symmetric
function, and let $\Gamma_k$ denote the connected component of 
$\{\lambda \in \Bbb R^n | \sigma_k (\lambda) > 0\}$ containing the
positive cone $\{\lambda \in \Bbb R^n|\lambda_1, \cdots, \lambda_n >
0\}$.  Then $(f, \Gamma) = (\sigma^{1/k}_k, \Gamma_k)$
satisfies (\ref{2})-(\ref{6}).

The special case of Problem \ref{MQuestClosed} 
for  
$(f, \Gamma)=(\sigma_1, \Gamma_1)$ is the 
Yamabe problem in the so-called positive case. The answer was proved affirmative through the works of Yamabe himself \cite{Yamabe60}, Trudinger \cite{Trudinger68}, Aubin \cite{Aubin76} and Schoen \cite{Schoen84}. Different solutions to the Yamabe problem in the case $n \leq 5$ and in the case $(M,g)$ is locally conformally flat were later given by Bahri and Brezis \cite{BB} and Bahri \cite{B}. In \cite{Schoen91}, Schoen proved a positive answer for Problem \ref{MQuest-extra} when $(M, g)$ is locally conformally flat and conjectured that the answer would also
be  positive for general Riemannian manifolds.
The conjecture was proved in dimensions $n\le 7$
by Li and Zhang  \cite{LiZhang05}
and Marques \cite{Marques05} independently.
For $n=3,4,5$, see works of
Li and Zhu \cite{LiZhu99},
Druet \cite{Druet03,Druet04}
and Li and Zhang \cite{LiZhang04}.
For $8\le n\le 24$,
it was proved that the answer to 
 Problem \ref{MQuest-extra} is positive provided 
 that the positive mass theorem holds in
these dimensions;
see Li and Zhang \cite{LiZhang05,LiZhang07}
for $8\le n\le 11$,
and  Khuri, Marques and Schoen
\cite{KMS09} for $12\le n\le 24$.  On the other hand,
the answer to  Problem \ref{MQuest-extra} is negative in
dimension $n\ge 25$; see
 Brendle \cite{Brendle08}
for $n\ge 52$, and
 Brendle  and Marques \cite{BrendleM09}
for $25\le n\le 51$.

Fully nonlinear elliptic equations involving
$f(\lambda(\nabla^2 u))$ were investigated in the
classic paper of Caffarelli, Nirenberg and Spruck \cite{C-N-S-Acta}.
Fully nonlinear elliptic equations involving the Schouten tensor and applications
to geometry and topology have been studied extensively in
and after the pioneering works of Viaclovsky \cite{Viac00-Duke, Viac00-TrAMS, Viac02} and Chang, Gursky and Yang
\cite{CGY02-AnnM, CGY02-JAM, CGY03-IHES, CGY03-IP}.
Extensions, as well as
developments of new methods, have been made by
Guan and Wang \cite{GW03-JRAM}, 
Li and Li \cite{LiLi03,LiLi05},
Gursky and Viaclovsky \cite{GV04, GV07},  
 Ge and Wang \cite{GeWang06},
 Sheng, Trudinger and Wang \cite{STW07},
 Trudinger and Wang \cite{TW09, TW10},
among others. Nevertheless, Problem \ref{MQuestClosed} is largely open for $2 < k < \frac{n}{2}$ and Problem \ref{MQuest-extra} is largely open for $2 \leq k < \frac{n}{2}$.

In \cite{LiNgBocher}, we began our study on Problem \ref{MQuest-extra}. In that paper, we restricted our attention to a locally conformally flat setting and established various asymptotic behavior near isolated singularities of the degenerate elliptic equation which arises naturally in the study of \eqref{9}, namely
\[ 
\lambda(A_{g_u}) \in \partial\Gamma \text{ in a punctured ball}. 
\] 
In this sequel to \cite{LiNgBocher}, we study compactness
of solutions of (\ref{9}). We consider
the following equation with a more general right hand side:
\begin{equation}
f(\lambda(A_{g_u})) =\psi(x), \quad \lambda (A_{g_u}) \in \Gamma,
\quad \mbox{ on } M,
\label{9'}
\end{equation}
where $\psi$ is a  given positive  smooth function  on $M$. We prove the following compactness result for \eqref{9'} under an additional assumption of a lower Ricci bound.

\begin{theorem}\label{Thm:Main} 
Let $(f, \Gamma)$ satisfy
(\ref{2})-(\ref{6}), $(M,g)$ be
a compact, smooth
Riemannian manifold of dimension $n\ge 3$, and $\psi$
be a positive smooth function on $M$. Assume that $(M,g)$ is not conformally equivalent to the standard sphere. For any $\alpha \geq 0$, either the set
\[
\mathfrak{S}_\alpha := \Big\{u \in C^2(M):  u \text{ is positive, satisfies \eqref{9'} and  $\Ric_{g_u} \geq - (n-1)\alpha^2\,g_u$}\Big\}
\]
is empty, or there exists $C = C(M,g,f,\Gamma,\psi,\alpha) > 0$ such that
\[
\|\ln u\|_{C^5(M)} \leq C \text{ for all } u \in \mathfrak{S}_\alpha.
\]
\end{theorem} 

\begin{remark}
Specific smoothness assumptions on $f$ and $\psi$ that ensure $C^{k,\beta}$ estimates for $u$ can be made precise from our proof. We however decided not to do so to keep the exposition clearer.
\end{remark}
 
Along the proof, the constant 
\begin{equation} 
\muGp \in [0,n-1] \text{ is the unique number such that } (-\muGp ,1, . . .1) \in \partial\Gamma, 
	\label{muGpDef} 
\end{equation} 
which was introduced in \cite{LiNgBocher}, plays an 
important role. $\muGp$ is well-defined thanks to \eqref{2} and \eqref{3}. For $\Gamma = \Gamma_k$, we have %
\[ 
\mu_{\Gamma_k}^+ = \frac{n-k}{k} \text{ for }1 \leq k \leq n. 
\] 
In particular, 
\[ 
\left\{\begin{array}{ll} 
\mu_{\Gamma_k}^+ > 1 &\text{ if } k < \frac{n}{2},\\ 
\mu_{\Gamma_k}^+ = 1 &\text{ if } k = \frac{n}{2},\\ 
\mu_{\Gamma_k}^+ < 1 &\text{ if } k > \frac{n}{2}. 
\end{array}\right. 
\] 
Also, we note that  
\begin{equation} 
\text{$(-\mu, 1, \ldots, 1) \in \Gamma$ for $\mu < \muGp$, and $(-\mu, 1, \ldots, 1) \in \RR^n \setminus \bar \Gamma$ for $\mu > \muGp$}. 
	\label{muGpProp} 
\end{equation} 
 
In our arguments, there is 
a key difference between the case $\muGp \leq 1$ and $\muGp > 1$. Note that for $\muGp \leq 1$, the lower Ricci bound assumption  with $\alpha = 0$ is 
satisfied automatically 
 (see \cite{G-V-W}). With the help of Theorem \ref{Thm:Main}, in fact its generalized version Theorem \ref{Thm:Cptnessmu<1} which allows estimates to hold uniformly along a homotopy connecting equation \eqref{9'} to a subcritical one, we obtain the following existence result for \eqref{9'} with $\muGp \leq 1$, which includes the $\sigma_k$-Yamabe problem for $k \geq \frac{n}{2}$.

\begin{theorem}\label{Thm:Existence} 
Let $(f, \Gamma)$ satisfy
(\ref{2})-(\ref{6}), $(M,g)$ be a compact, smooth
Riemannian manifold of dimension $n\ge 3$
satisfying $\lambda(A_g)\in \Gamma$ on $M$, and $\psi$ be a positive smooth function on $M$. If $\muGp \leq 1$, then there exists a 
smooth positive solution $u$ of \eqref{9'}.
Moreover, if $(M,g)$ is not conformally equivalent to the standard sphere, then all solutions $u$  of \eqref{9'}
satisfy $\|\ln u\|_{ C^5(M,g) }\le C$ for some constant
$C$ depending only on $(f,\Gamma)$, $(M,g)$ and $\psi$.
\end{theorem}  

\begin{remark}  In fact, the degree of all solutions in
the above theorem is equal to $-1$,
as proved in Section \ref{degree}.
\end{remark}
  
Prior to our work, when $(f,\Gamma) \neq (\sigma_1, \Gamma_1)$,
the state of the art was, roughly speaking, as follows:
\begin{enumerate}[(i)] 
  \item  A positive  answer to 
both Problem \ref{MQuestClosed} and Problem \ref{MQuest-extra} for $(f,\Gamma) = (\sigma_2^{1/2},\Gamma_2)$ in dimension $n = 4$, see \cite{CGY02-JAM}.
\item A positive  answer to 
both Problem \ref{MQuestClosed} and Problem \ref{MQuest-extra} when
 $(f, \Gamma)$ satisfies (\ref{2})-(\ref{6}) and 
$(M^n, g)$ is locally conformally flat, see
 \cite{GW03-JRAM, LiLi03,  LiLi05},
 
\item A positive  answer to Problem \ref{MQuestClosed} for
$(f,\Gamma)=(\sigma_2^{\frac 12}, \Gamma_2)$, see \cite{GeWang06, STW07}.

\item A positive answer to 
both Problem \ref{MQuestClosed} and Problem \ref{MQuest-extra} for
 $(f,\Gamma)=(\sigma_k^{\frac 1k}, \Gamma_k)$ with $
k > \frac n2$, see \cite{GV04, GV07, TW09}. 
 
\item A positive answer to 
both Problem \ref{MQuestClosed} and Problem \ref{MQuest-extra} for $(f,\Gamma)=(\sigma_k^{\frac 1k}, \Gamma_k)$ and $k = \frac n2$ was given in \cite{TW10}, though we do not follow the proof.

\end{enumerate} 
 
On one hand, our results cover the above statements for $k \geq \frac{n}{2}$. On the other hand, they give new results for the classical Yamabe problem; we note that the answer to
 Problem \ref{MQuest-extra} when $(f,\Gamma) = (\sigma_1,\Gamma_1)$ is negative for $n\ge 25$.

In our proof of
Theorem \ref{Thm:Main}, we use the Ricci lower bound to show directly that there cannot be more than one blow-up point and, if there is a blow-up point, there is no bubble accumulation. This is very different from existing compactness arguments in the literature. 
For this step in the case $\mu^+_\Gamma\le 1$, we make use
of Bishop's comparison theorem, 
a Riemannian version of Hawking's singularity theorem in
relativity, a Liouville theorem in \cite{LiLi05} and local gradient estimates for solutions of \eqref{9'} (\cite{GW03-IMRN, Chen05, Li09-CPAM, Wang06}).
In the case $\mu^+_\Gamma>1$, 
we also make use of an isoperimetric inequality of
B\'erard, Besson and Gallot in \cite{BBG85-InvM}.
The rest of the argument is to obtain estimates for the blow-up limit so that one can apply the Bishop-Gromov comparison theorem as in \cite{GV07} to reach the conclusion. This part makes use of a symmetry result in \cite{Li09-CPAM} for solutions
of $\lambda(A_{u^{\frac{4}{n-2}}g_{\rm flat}})\in \partial \Gamma$ on a punctured space and certain constructions
 of sub-solutions and super-solutions of 
the equation $\lambda(A_{g_u})\in
\partial \Gamma$.
Much of the constructions
is   based on knowledge
obtained
 in our earlier work \cite{LiNgBocher} in the locally conformally flat case.

We also give a second proof of Theorem \ref{Thm:Main} which avoids the use of the Riemannian version of Hawking's singularity theorem, though this proof is somewhat more elaborate.

The rest of the paper is structured as follows. In Section \ref{Sec:Hawking} we start with a Riemannian version of Hawking's singularity theorem in relativity. In Section \ref{Sec:Main}, we give the proof of Theorem \ref{Thm:Main}. In Section \ref{general}, we present a generalization of Theorem \ref{Thm:Main} which includes subcritical equations. In Section \ref{degree}, we prove Theorem \ref{Thm:Existence}. In Section \ref{secondproof} we give a second proof of Theorem \ref{Thm:Main}. In the appendix, we prove some auxiliary results that were needed in the body of the paper. 
 
\section{Preliminary}\label{Sec:Hawking} 
 
We start with a Riemannian version of Hawking's singularity theorem (see e.g. \cite[page 271]{HawkingEllis}): 
 
\begin{proposition}\label{Prop:Hawking} 
Let $(N^n,g)$ be a complete smooth Riemannian manifold with smooth boundary $\partial N$. If $\Ric_g \geq - (n-1)\alpha^2$ for some $\alpha \geq 0$ and if the mean curvature $H$ of $\partial N$ with respect to its inward pointing normal satisfies $H > (n - 1)c_0 > (n - 1)\alpha$, then  
\[ 
d_g(x,\partial N) \leq U(\alpha,c_0) \text{ for all } x \in N, 
\] 
where $d_g$ denotes the distance function induced by $g$ and 
\[ 
U(\alpha,c_0) = \left\{\begin{array}{ll} 
	\frac{1}{c_0} &\text{ if } \alpha = 0,\\ 
	\frac{1}{\alpha}\,\coth^{-1}\big(\frac{c_0}{\alpha}\big) & \text{ if } \alpha > 0. 
\end{array}\right. 
\] 
\end{proposition} 
 
\bproof The proof is a standard argument using the second 
variation formula of arc-length. We include it here for completeness. 
 
Fix $x \in N$ and let $\gamma(t): [0,a] \rightarrow N$ be a unit speed geodesic with $\gamma(0) \in \partial N$ and $\gamma(a) = x$ such that $a = d_g(x,\partial N)$. We need to show that $a \leq U(\alpha,c_0)$. 
 
Clearly $\gamma'(0) \perp T_{\gamma(0)} (\partial N)$. Choose an orthonormal frame $E_1, \ldots, E_n$ at $\gamma(0)$ such that $E_1 = \gamma'(0)$ and extend it along $\gamma$ by parallel transport.  
 
Fix some $i \in \{2, \ldots, n\}$ and let $\gamma_i: [-\delta, \delta] \times [0,a] \rightarrow N$ be a variation of $\gamma$ such that $\gamma_i(s,0) \in \partial N$, $\gamma_i(s,a) = x$, $\gamma_i(0,t) = \gamma(t)$ and 
\[ 
\frac{\partial}{\partial s}\Big|_{s = 0} \gamma_i(s,t) = V_i = f(t)\,E_i(t), 
\] 
for some $f(t)$ which will be specified later. By the second variation formula of arc-length (see e.g. \cite[Theorem 2.5]{ChavelModernIntro}) and the minimizing property of $\gamma$, we have 
\begin{align*} 
0  
	&\leq \frac{1}{2}\frac{d^2}{ds^2}\Big|_{s = 0} \mathrm{Length}(\gamma_i(s,\cdot)) \\ 
	&= -g(\nabla_{V_i}{V_i},\gamma'(0)) + \int_0^a \Big[g\Big(\frac{D}{dt}V_i, \frac{D}{dt}V_i\Big) - g(V_i,R(\gamma', V_i)\gamma') \Big]\,dt\\ 
	&= -f(0)^2\,g(\nabla_{E_i}{E_i},\gamma'(0))  + \int_0^a \Big[|f'(t)|^2 - f(t)^2\,g(E_i, R(\gamma', E_i)\gamma')\Big]\,dt. 
\end{align*} 
Summing over $i$ and using our hypotheses on the Ricci curvature and the mean curvature, we obtain 
\begin{align*} 
0  
	&\leq -f(0)^2\,H(\gamma(0)) + \int_0^a \Big[(n-1)|f'(t)|^2 - f(t)^2\,\Ric(\gamma',\gamma')\Big]\,dt\\ 
	&\leq -(n - 1)\,f(0)^2\,c_0 + (n-1)\int_0^a \Big[|f'(t)|^2 + \alpha^2\,f(t)^2\Big]\,dt. 
\end{align*} 
Optimizing the right hand side subjected to $f(0) = 1$ and $f(a) = 0$ leads to 
\[ 
f(t) = \left\{ \begin{array}{ll} 
	1 - \frac{t}{a} & \text{ if } \alpha = 0,\\ 
	-\frac{\sinh(\alpha(t - a))}{\sinh(\alpha a)} & \text{ if } \alpha > 0. 
\end{array}\right. 
\] 
Using this choice of $f$, we arrive at 
\[ 
0 \leq \left\{ \begin{array}{ll} 
	-c_0 + \frac{1}{a} & \text{ if } \alpha = 0,\\ 
	- c_0 + \alpha \coth(\alpha a) & \text{ if } \alpha > 0. 
\end{array}\right. 
\] 
The conclusion follows. 
\eproof

\section{Proof of Theorem \ref{Thm:Main}}\label{Sec:Main}

Without loss of generality, we assume that $f(\lambda(A_{g_{\rm can}})) = 1$ on $\SSphere^n$, where $g_{\rm can}$ is the standard metric on $\SSphere^n$. 
 
For simplicity, we present the proof for $\psi \equiv 1$. The proof for general $\psi$ requires only minor modifications. 

Assume that $\mathfrak{S}_\alpha$ is non-empty for some $\alpha \geq 0$. Assume for the moment that we have established:
\begin{equation} 
\max_M u \leq C \text{ for all } u \in \mathfrak{S}_\alpha.
\label{maxestimate}
\end{equation}
By \cite[Theorem 1.10]{Li09-CPAM} and \cite[Theorem 1.20]{LiLi03}, this implies that
\begin{equation}
 \max_M [|\nabla_g\ln u|+|\nabla^2 \ln u|]\le C \text{ for all } u \in \mathfrak{S}_\alpha.
	\label{maxestCons}
\end{equation}
The desired estimate on $\ln u$ follows,
in view of Evan-Krylov's and Schauder's estimates, once we can show that
\begin{equation}
\min_M u \geq \frac{1}{C} \text{ for all } u \in \mathfrak{S}_\alpha.
	\label{minestimate}
\end{equation}

To prove \eqref{minestimate} assume by contradiction that there is a sequence $u_i$ in $\mathfrak{S}_\alpha$ such that 
\begin{equation}
\min_M u_i \rightarrow 0.
	\label{Blowdown}
\end{equation}
By definition, the metrics $g_i := u_i^{\frac{4}{n - 2}}\,g$ satisfy 
\begin{equation} 
f(\lambda(A_{g_i})) = 1, \qquad \lambda(A_{g_i}) \in \Gamma \;\text{ on } M. \label{9i} 
\end{equation} 
Noting that $g_i = (\frac{u_i}{u_1})^{\frac{4}{n-2}}\,g_1$ and evaluating (\ref{9i}) at a maximum point
$\bar x_i$  of
$\frac{u_i}{u_1}$, we obtain 
\[
1 \ge \Big(\frac{u_i(\bar x_i)}{u_1(\bar x_i)}\Big)^{ -\frac 4{n-2}}
 f(\lambda(A_{g_1}(\bar x))
	= \Big(\frac{u_i(\bar x_i)}{u_1(\bar x_i)}\Big)^{ -\frac 4{n-2}},
\]
which implies
$\max_M u_i \ge \min_M u_1$. Recalling \eqref{maxestCons}, we see that $\min_M u_i \geq \frac{1}{C} \min_M u_1$, which contradicts \eqref{Blowdown}.

We turn to the proof of \eqref{maxestimate}. Arguing by contradiction, assume that there is a sequence $u_i$ of smooth positive functions on $M$ such that the metrics $g_i = u_i^{\frac{4}{n - 2}}\,g$ satisfy \eqref{9i} and 
\begin{equation} 
\Ric_{g_i} \geq -(n-1)\alpha^2\,g_i \;\text{ on M}, \label{RicLB} 
\end{equation} 
but 
\begin{equation} 
M_i := u_i(x_i) = \max_M u_i \rightarrow \infty. 
	\label{Blowup} 
\end{equation} 
We can assume that $x_i \rightarrow x_\infty$ in the topology induced by $g$. We would like to arrive at a contradiction. 
 
The proof will be divided according to whether $\muGp \leq 1$ or $\muGp > 1$. In each case, the proof consists of six steps.

\subsection{The case $\muGp \leq 1$}\label{musmall}

\noindent \underline{Step 1:} We show that $x_\infty$ is the unique blow-up point of $\{u_i\}$. In fact, we show a stronger result:
For some constant $C$ independent of $i$, 
\begin{equation} 
u_i(x)  \leq Cd_g(x,x_i)^{-\frac{n-2}{2}} \text{ for all } x \in M \setminus \{x_i\}.
	\label{IsolatedBlowupPt} 
\end{equation} 

The key ingredient for this step is the following result.

\begin{lemma}\label{Lem:NearBlowupPt} 
Assume for some $C_1 \geq 1$, $K_i \rightarrow \infty$ and $y_i \in M$ that 
\begin{equation} 
u_i(y_i) \rightarrow \infty \text{ and } \sup \Big\{u_i(y): d_g(y,y_i) \leq K_i\,u_i(y_i)^{-\frac{2}{n-2}}\Big\} \leq C_1\,u_i(y_i) .
	\label{NBP:Cond} 
\end{equation} 
Then for any $0 < \mu < 1$, there exists $K = K(C_1,\mu)$ such that for all sufficiently large $i$
\[ 
\textrm{Vol}_{g_i}\Big(\Big\{y: d_g(y,y_i) \leq Ku_i(y_i)^{-\frac{2}{n-2}}\Big\}\Big) \geq (1 - \mu) \textrm{Vol}_{g_i}(M). 
\] 
\end{lemma}

\bproof Define, for $p \in \RR^n$ and $a > 0$, 
\begin{equation}
U_{a,p}(x) = c\,\Big(\frac{a}{1 + a^2|x - p|^2}\Big)^{\frac{n-2}{2}}, \qquad x \in \RR^n,  
	\label{UapDef}
\end{equation}
where $c = 2^{\frac{n-2}{2}}$. 
 
Write 
\[ 
\SSphere^n = \{z = (z_1, \ldots, z_{n+1}) \in \RR^{n+1} | z_1^2 + \ldots z_{n+1}^2 = 1\}. 
\] 
Let $(x_1, \ldots, x_n) \in \RR^n$ be the stereographic projection coordinates of $z \in \SSphere^n$, i.e. 
\[ 
z_i = \frac{2x_i}{1 + |x|^2} \text{ for }  1 \leq i \leq n, \text{ and }z_{n+1} = \frac{|x|^2 - 1}{|x|^2 + 1}. 
\] 
Then 
\[ 
g_{\rm can} = |dz|^2 = \Big(\frac{2}{1 + |x|^2}\Big)^2\,|dx|^2 = U_{1,0}^{\frac{4}{n-2}}\,|dx|^2. 
\] 
Thus, by conformal invariance, we have 
\[ 
f(\lambda(A_{U_{a,p}^{\frac{4}{n-2}}\,g_{\rm flat}})) = 1 \text{ on } \RR^n, 
\] 
where $g_{\rm flat} = |dx|^2$ is the standard Euclidean metric on $\RR^n$.

Define a map $\Phi_i: \RR^n \approx T_{y_i}(M,g) \rightarrow M$ by 
\[ 
\Phi_i(x) = \exp_{y_i} \frac{c^{\frac{2}{n-2}}\,x}{u_i(y_i)^{\frac{2}{n-2}}},
\] 
and let  
\[ 
\tilde u_i(x) = \frac{c}{u_i(y_i)}u_i \circ \Phi_i(x), \qquad x  \in \RR^n. 
\] 
Then $\tilde u_i$ satisfies 
\[ 
f(\lambda(A_{\tilde u_i^{\frac{4}{n-2}}\,\tilde h_i})) = 1 \text{ and } \lambda(A_{\tilde u_i^{\frac{4}{n-2}}\,\tilde h_i})) \in \Gamma \text{ on } \{|x| < \delta_0\,c^{-\frac{2}{n-2}}\,u_i(y_i)^{\frac{2}{n-2}}\}, 
\] 
where $\tilde h_i := \Big(\frac{c}{u_i(y_i)}\Big)^{-\frac{4}{n-2}}\Phi_i^* g$ and $\delta_0$ is the injectivity radius of $(M,g)$. It is clear that $\tilde h_i \rightarrow g_{\rm flat}$ on $C^3_{\rm loc}(\RR^n)$. Furthermore $\tilde u_i(0) = c$ and, by \eqref{NBP:Cond}, $\tilde u_i \leq C_1\,c$. By known $C^1$, $C^2$ estimates (\cite[Theorem 1.10]{Li09-CPAM} and \cite[Theorem 1.20]{LiLi03}), $\ln \tilde u_i$ is uniformly bounded in $C^2$ on any compact subset of $\RR^n$. By Evans-Krylov's theorem \cite{Evans,Krylov} and the Schauder theory, $\tilde u_i$ is uniformly bounded in $C^3$ on any compact subset of $\RR^n$and subconverges in $C^{2,\alpha}_{\rm loc}(\RR^n)$ to some positive $\tilde u_* \in C^2(\RR^n)$ which satisfies 
\[ 
f(\lambda(A_{\tilde u_*^{\frac{4}{n-2}}\,g_{\rm flat}})) = 1 \text{ and } \lambda(A_{\tilde u_*^{\frac{4}{n-2}}\,g_{\rm flat}})) \in \Gamma \text{ on } \RR^n.
\] 
By the Liouville theorem \cite[Theorem 1.3]{LiLi05}, we have $\tilde u_* = U_{a_*,x_*}$  
for some $a_* > 0$ and $x_* \in \RR^n$. Since $\tilde u_*(0) = \lim \tilde u_i(0) = c$ and $\tilde u_* \leq C_1\,c$, we have,
for some constant $C$ depending only on
$C_1$ and $n$, that  
\begin{equation} 
|x_*| \leq C \text{ and } C^{-1} \leq a_* \leq C. 
\end{equation} 
In particular, for any $R > 0$ and $\mu > 0$,  
\begin{equation} 
\|\tilde u_i - \tilde u_*\|_{C^2(\bar B_R)} \leq \mu \text{ for all sufficiently large $i$}. 
	\label{Eq:Liouville} 
\end{equation} 

It follows that the metrics $\tilde u_i^{\frac{4}{n-2}}\,\tilde h_i$ converge on compact subsets to the metric $\tilde u_*^{\frac{4}{n-2}}g_{\rm flat}$. Since $(B(0,r), \tilde u_i^{\frac{4}{n-2}}\,\tilde h_i)$ is isometric to $(\Phi_i(B(0,r)),g_i)$ for any $r > 0$, we obtain: For any $\epsilon > 0$, there exists $R = R(\epsilon,C_1) > 0$ such that  
\begin{enumerate}[(i)] 
\item $|\textrm{Vol}_{g_i}(\Phi_i(B(0,R))) - \textrm{Vol}(\SSphere^n, g_{\rm standard})| \leq C\,\epsilon^n$ for some $C$ independent of $i$ and $\epsilon$, 
\item and the mean curvature of the hypersurface $\partial \Phi_i(B(0,R))$ with respect to $g_i$ and the unit normal pointing away from $\Phi_i(B(0,R))$ is no smaller than $\frac{1}{\epsilon}$. 
\end{enumerate} 
Using (ii) and Proposition \ref{Prop:Hawking}, we see that 
\[ 
\textrm{diam}_{g_i}(M\setminus\Phi_i(B(0,R))) \leq C\,\epsilon. 
\] 
In view of \eqref{RicLB} and Bishop's theorem (see e.g. \cite[Theorem 3.9]{ChavelModernIntro}), this implies that 
\[ 
\textrm{Vol}_{g_i}(M\setminus\Phi_i(B(0,R))) \leq C\,\epsilon^n. 
\] 
Lemma \ref{Lem:NearBlowupPt} is established. 
\eproof

We are now in position to prove \eqref{IsolatedBlowupPt}. Assume that \eqref{IsolatedBlowupPt} is incorrect, then, for some $\tilde x_i \in M$, 
\begin{equation} 
u_i(\tilde x_i)\,d_g(x_i,\tilde x_i)^{\frac{n-2}{2}} = \max_M u_i\,d_g(x_i,\cdot)^{\frac{n-2}{2}} \rightarrow \infty. 
	\label{NegIBP} 
\end{equation} 
Since $(M,g)$ is compact, this implies that $u_i(\tilde x_i) \rightarrow \infty$.

Applying Lemma \ref{Lem:NearBlowupPt} to $C_1 = 1$, $y_i = x_i$ and $K_i = \delta\,u_i(x_i)^{\frac{2}{n-2}}$ with some small $\delta = \delta(M,g)$, we find 
\[ 
\textrm{Vol}_{g_i}\Big(\Big\{y: d_g(y,x_i) \leq Ku_i(x_i)^{-\frac{2}{n-2}}\Big\}\Big) \geq \frac{3}{4} \textrm{Vol}_{g_i}(M), 
\] 
where $K$ is some universal constant. Also, applying Lemma \ref{Lem:NearBlowupPt} again to $C_1 = 2^{\frac{n-2}{2}}$, $y_i = \tilde x_i$ and $K_i = \frac{1}{2}\,d(x_i,\tilde x_i)\,u_i(\tilde x_i)^{\frac{2}{n-2}}$, we obtain 
\[ 
\textrm{Vol}_{g_i}\Big(\Big\{y: d_g(y,\tilde x_i) \leq \tilde K\,u_i(\tilde x_i)^{-\frac{2}{n-2}}\Big\}\Big) \geq \frac{3}{4} \textrm{Vol}_{g_i}(M), 
\] 
where $\tilde K$ is another universal constant. On the other hand, since $u_i(x_i) \geq u_i(\tilde x_i)$, \eqref{NegIBP} implies that the sets  
\[ 
\Big\{y: d_g(y,x_i) \leq Ku_i(x_i)^{-\frac{2}{n-2}}\Big\} \text{ and }\Big\{y: d_g(y,\tilde x_i) \leq \tilde K\,u_i(\tilde x_i)^{-\frac{2}{n-2}}\Big\} 
\] 
are disjoint for all sufficiently large $i$. The last three conclusions are incompatible thus yield a contradiction. We have proved \eqref{IsolatedBlowupPt}. 
 
\begin{remark} 
The above argument also shows that 
\begin{align} 
\textrm{diam}_{g_i}(M) \rightarrow \textrm{diam}(\SSphere^n,g_{\rm can}), \label{diamBnd}\\ 
\textrm{Vol}_{g_i}(M) \rightarrow \textrm{Vol}(\SSphere^n,g_{\rm can}). \label{VolBnd} 
\end{align} 
\end{remark} 
 
\medskip 
\noindent\underline{Step 2:} We prove that 
\begin{equation} 
|\nabla_g^k \ln u_i(x)| \leq C\,d_g(x,x_i)^{-k} \text{ for } x \neq x_i, k = 1,2. 
	\label{SC2Est} 
\end{equation} 

By Step 1 and known estimates, we have, for any compact $K \subset M\setminus \{x_\infty\}$, there exists $N(K)$ such that 
\[ 
\|u_i\|_{C^2(K)} \leq C(K) \text{ for all } i \geq N(K). 
\] 
To prove the more precise form \eqref{SC2Est}, fix some $x \neq x_i$ and let $d_i = d_g(x,x_i)$. Define $\Psi_i: B_{1/2} \subset \RR^n \approx T_x(M,g) \rightarrow M$ by 
\[ 
\Psi_i(y) = \exp_x \Big(d_i\,y\Big). 
\] 
Then the metrics $d_i^{-2}\,\Psi_i^* g$ on $B_1$ have a uniform injectivity radius and curvature bound. Furthermore, the function  
\begin{equation}
w_i(y) := d_i^{\frac{n-2}{2}}\,u_i(\Psi_i(y))  
	\label{Eq:wiDef}
\end{equation}
is uniformly bounded by \eqref{IsolatedBlowupPt} and satisfies 
\[ 
f(\lambda(A_{w_i^{\frac{4}{n-2}}\,d_i^{-2}\,\Psi_i^* g})) = 1, \lambda(A_{w_i^{\frac{4}{n-2}}\,d_i^{-2}\,\Psi_i^* g}) \in \Gamma \text{ in } B_{1/2}. 
\] 
Thus, by \cite[Theorem 1.10]{Li09-CPAM} and \cite[Theorem 1.20]{LiLi03},  
\[ 
|\nabla_{d_i^{-2}\,\Psi_i^* g}^k \ln w_i| \leq C \text{ in } B_{1/4}, k = 1,2. 
\] 
Returning to $u_i$, we obtain \eqref{SC2Est}. 
 
\medskip 
\noindent\underline{Step 3:} We show that $u_i$ converges in $C^0_{\rm loc}(M\setminus\{x_\infty\})$ to $0$. 
 
By Step 1 and Step 2, $u_i$ converges uniformly away from $x_\infty$ to some non-negative limit $u_\infty$. Lemma \ref{Lem:NearBlowupPt} and the volume bound \eqref{VolBnd} shows that, for any $\delta > 0$, 
\[ 
\int_{\{x: d_g(x,x_i) > \delta\}} u_i^{\frac{2n}{n-2}}\,dv_g = \textrm{Vol}_{g_i}\Big(\{x: d_g(x,x_i) > \delta\}\Big) \rightarrow 0. 
\] 
Sending $i$ to infinity, we see that $u_\infty \equiv 0$. 
 
Since $u_\infty \equiv 0$, in order to obtain a useful blow-up limit, we need to rescale the sequence $\{u_i\}$. Fix some point $p \in M \setminus \{ x_\infty\}$ and let  
\[ 
v_i(x) = u_i(p)^{-1}\,u_i(x). 
\] 
Note that $v_i(p) = 1$. Thus, by Step 2, $v_i$ subconverges, for every $0 < \alpha < 1$, in $C^{1,\alpha}_{\rm loc}(M \setminus \{x_\infty\},g)$ to some positive function $v_\infty \in C^{1,1}_{\rm loc}(M \setminus \{x_\infty\},g)$, which satisfies $v_\infty(p) = 1$ and  
\begin{equation} 
|\nabla_g^k \ln v_\infty(x)| \leq C\,d_g(x,x_\infty)^{-k} \text{ in } M \setminus \{x_\infty\}, k = 1, 2. 
	\label{vIC2Est} 
\end{equation} 
Furthermore, as $u_i(p) \rightarrow 0$, $v_\infty$ satisfies, in view of the equation satisfied by $u_i$, in viscosity sense 
\begin{equation} 
\lambda(A_{g_{v_\infty}}) \in \partial \Gamma \text{ in } M\setminus \{x_\infty\}. 
	\label{vIDegEqn} 
\end{equation} 

\medskip 
\noindent\underline{Step 4:} We show that, for some constant $a \geq 0$,
\begin{equation} 
\lim_{x \rightarrow x_\infty} v_\infty(x)\,d_g(x,x_\infty)^{n-2} = a \in [0,\infty). 
	\label{vIAsympPrecise} 
\end{equation} 

In the case where the background is locally conformally flat, estimate \eqref{vIAsympPrecise} was derived in \cite[Theorem 1.10]{LiNgBocher}. We adapt the argument therein to the case at hand. 
 
We first show that 
\begin{equation} 
a:= \lim_{r \rightarrow 0} r^{n-2} \min_{\partial B_g (x_\infty,r)} v_\infty  \text{ is  finite}. 
	\label{vILiminf} 
\end{equation} 
In view of the gradient estimate in \eqref{vIC2Est}, we deduce from the above that 
\begin{equation} 
A:=\limsup_{x \rightarrow x_\infty} v_\infty(x)\,d_g(x,x_\infty)^{n-2} \text{ is  finite}. 
	\label{vILimsup} 
\end{equation} 
Estimate \eqref{vILiminf} is a consequence of the super-harmonicity of $v_\infty$ (with respect to the conformal Laplacian); see the lemma below.

\begin{lemma}\label{Lem:SupH} 
Let $\Omega$ be an open neighborhood of a point $p \in M$. Let $L_g = \Delta_g - c(n)\,R_g$ denote the conformal Laplacian of $g$. If $w \in LSC(\Omega \setminus \{p\})$ is a non-negative function in $\Omega \setminus \{p\}$ and satisfies $L_gw \leq 0$ in the viscosity sense in $\Omega \setminus \{0\}$, then  
\[ 
\lim_{r \rightarrow 0} r^{n-2}\,\min_{\partial B_g(p,r)} w \text{ exists and is finite}. 
\] 
\end{lemma}

\bproof When $L_g$ is the Laplacian on a Euclidean space, the lemma is classical.  
 
Using normal coordinates, we can assume that $p = 0$, $\Omega = B_1$, $d_g(x,p) = |x| + O(|x|^2)$ and 
\[ 
L_g = a_{ij}(x)\,D_{ij} + b_i(x)\,D_i + c(x) 
\] 
where $a_{ij}(x) = \delta_{ij} + O(r), b_i(x) = O(1), c(x) = O(1)$, the big `O' notation is meant for $x$ close to the origin, and $r = |x|$.  
 
Set $\underline{w}(r) = \min_{\partial B_r} w$.

A calculation gives 
\begin{align*} 
L_g r^{2 - n} 
	&= O(r^{1-n}),\\ 
L_g r^{\frac{5}{2} - n} 
	&= -\frac{1}{2}(n - \frac{5}{2})r^{\frac{1}{2} - n} + O(r^{\frac{3}{2} - n}). 
\end{align*} 
Thus, for some $K$ sufficiently large and $\delta$ sufficiently small, the function $G(x) = r^{2 - n} - Kr^{\frac{5}{2} - n} - (\delta^{2 - n} - K\delta^{\frac{5}{2} - n})$ is non-negative and  
\[ 
L_g G \geq 0 \text{ in } B_\delta \setminus \{0\}. 
\] 

By the maximum principle (which holds for $\delta$ sufficiently small), we have for any $0 < \rho < \delta$, 
\[ 
w \geq \frac{\underline{w}(\rho)}{G(\rho)}G \text{ in } \{\rho \leq r \leq \delta\}, 
\] 
and so
\[ 
\underline{w}(r) \geq \frac{\underline{w}(\rho)}{G(\rho)}G(r) \text{ for all } r \in (\rho,\delta). 
\] 
It follows that the function $G^{-1}\,\underline{w}$ is increasing in $(0,\delta)$, in particular 
\[ 
\lim_{r \rightarrow 0} \frac{\underline{w}(r)}{G(r)} \text{exists and is finite}. 
\] 
The conclusion follows. 
\eproof

We proceed with the proof of \eqref{vIAsympPrecise}. It remains to show that $A = a$. Arguing indirectly, assume that $A > a$. Then we can find a sequence $z_i \rightarrow x_\infty$ such that, for some $\epsilon > 0$, 
\[ 
A + \epsilon \geq d_g(z_i,x_\infty)^{n - 2}\,v_\infty(z_i) \geq a + 2\epsilon 
\] 
while, in view of \eqref{vILiminf}, 
\[ 
d_g(z_i,x_\infty)^{n - 2}\,\min_{d_g(z,x_\infty) = d_g(z_i,x_\infty)} v_\infty(z) \leq a + \epsilon. 
\] 

Let $R_i = d_g(z_i,x_\infty)^{-1}$ and some sufficiently small $\delta > 0$. Define $\Xi_i: B_{\delta R_i} \subset \RR^n \approx T_{x_\infty}(M,g) \rightarrow M$ by 
\[ 
\Xi_i(y) = \exp_{x_\infty} (R_i^{-1}\,y). 
\] 
As before, $R_i^2\,\Xi_i^* g$ converges on compact subsets to the Euclidean metric $g_{\rm flat}$. 
 
Set 
\[ 
\hat v_i(y) = R_i^{2 - n} v_\infty \circ \Xi_i(y). 
\] 
Then $\hat v_i \in C^{1,1}_{\rm loc}(B_{\delta R_i} \setminus \{0\})$ satisfies 
\begin{equation} 
\left\{\begin{array}{l}\lambda\Big(A_{\hat v_i^{\frac{4}{n-2}}\,R_i^2\,\Xi_i^* g}\Big) \in \partial \Gamma \text{ in } B_{\delta\,R_i}\setminus \{0\},\\ 
\min_{\partial B_1} \hat v_i \leq a + \epsilon \text{ and } \max_{\partial B_1} \hat v_i \geq a + 2\epsilon. 
\end{array}\right. 
	\label{EqX1} 
\end{equation} 
Also, by \eqref{vIC2Est}, we have for some constant $C$ independent of $i$, 
\[ 
|\nabla_{R_i^2\,\Xi_i^* g}^k \ln \hat v_i(x)| \leq C\,|x|^{-k} \text{ in } B_{\delta\,R_i} \setminus \{0\}, k = 1, 2. 
\] 
Using the above estimate and the second line of \eqref{EqX1}, we see that, up to a subsequence, $\hat v_i$ converges uniformly on compact subsets of $\RR^n \setminus \{0\}$ to some limit $\hat v_* \in C^{1,1}_{\rm loc}(\RR^n \setminus \{0\})$, which, by the first line of \eqref{EqX1}, satisfies in the viscosity sense 
\begin{equation} 
\lambda\Big(A_{\hat v_*^{\frac{4}{n-2}}\,g_{\rm flat}}\Big) \in \partial \Gamma \text{ in } \RR^n\setminus \{0\}. 
	\label{LimitDegEqn} 
\end{equation} 
By \cite[Theorem 1.18]{Li09-CPAM}, $\hat v_*$ is radially symmetric, i.e. $\hat v_*(y) = \hat v_*(|y|)$. This results in a contradiction with the second line of \eqref{EqX1} and the convergence of $\hat v_i$ to $\hat v_*$. We have thus established \eqref{vIAsympPrecise}.

\medskip 
\noindent\underline{Step 5:} We now show that $a > 0$. 
 
So far, we have not used the assumption that $\muGp \leq 1$. We will show that there exists some $C > 0$ and $r_1 > 0$ independent of $i$ such that, for any small $\delta > 0$, there holds for some $K = K(\delta) > 0$ and $N = N(\delta) > 0$ that
\begin{equation} 
v_i(x) \geq \frac{1}{C\,d_g(x,x_i)^{n - 2 - 2\delta}} \text{ in } \{K\,u_i(x_i)^{-\frac{2}{n-2}} \leq d_g(x,x_i) \leq r_1\} \text{ for all } i \geq N. 
	\label{Eq:mu=1viLB} 
\end{equation} 
Clearly this implies that $v_\infty \geq \frac{1}{C\,d_g(x,x_\infty)^{n-2-2\delta}}$ in $\{0 < d_g(x,x_\infty) \leq r_1\}$ for all sufficiently small $\delta >0$, which in turns implies that $a > 0$. 

In normal coordinates at $x_i$, let $r = |x|$. We will use the following lemma, whose proof can be found in Appendix \ref{App:SS}. 
 
\begin{lemma}\label{Lem:mu=1SubSol} 
Assume that $\Gamma$ satisfies \eqref{2}, \eqref{3} and that $\muGp \leq 1$.
There exists some small $r_1>0$ depending
only on $(M,g)$ such that for all $0<\delta
<\frac 14$,
 the function $\underline{v}_\delta := 
r^{-(n-2-2\delta)}e^r$ satisfies
\[ 
\lambda(A_{(\underline{v}_\delta)^{\frac{4}{n-2}}g}) \in \RR^n \setminus \bar\Gamma \text{ in } \{0 < r < r_1\}. 
\] 
\end{lemma} 
 
From $v_i(p) = 1$ and \eqref{SC2Est}, there exists some positive constant $C$ independent of $i$ and $\delta$ such that 
\[ 
v_i \geq \frac{1}{C}\underline{v}_\delta \text{ on } \{r = r_1\}. 
\] 
For some $K = K(\delta) > 0$ to be fixed, let 
\[ 
\bar\beta = \sup\Big\{0 < \beta < \frac{1}{C}: v_i \geq \beta \underline{v}_\delta \text{ in } \{r_i := K\,u_i(x_i)^{-\frac{2}{n-2}} < r < r_1\}\Big\}. 
\] 
To prove \eqref{Eq:mu=1viLB}, it suffices to show that $\bar\beta = \frac{1}{C}$. Arguing indirectly, assume that $\bar\beta < \frac{1}{C}$. Then, in view of the equation satisfied by $v_i$, Lemma \ref{Lem:mu=1SubSol} and the comparison principle, there exists $\hat x_i$ with $|\hat x_i| = r_i$ such that  
\[ 
v_i(\hat x_i) = \bar\beta \underline{v}_\delta(\hat x_i). 
\] 
Since $v_i \geq \bar\beta \underline{v}_\delta$ in $\{r_i < r < r_1\}$, it follows that 
\begin{equation} 
\partial_r \ln v_i(\hat x_i) \geq \partial_r \ln(\bar\beta\underline{v}_\delta)(\hat x_i) = \partial_r \ln \underline{v}_\delta(\hat x_i) > - \frac{n - 2 - \delta}{r_i}. 
	\label{Eq:mu=1SSCons} 
\end{equation} 

On the other hand, by \eqref{Eq:Liouville} with $y_i = x_i$ (note that $\tilde u_* = U_{a,p}$ with $a = 1$ and $p = 0$ in the present case), 
\[ 
\Big|\frac{c^{\frac{2}{n-2}}}{u_i(x_i)^{\frac{2}{n-2}}}\partial_r \ln v_i(\hat x_i) + \frac{(n-2)c^{-\frac{2}{n-2}}\,K}{1 + c^{-\frac{4}{n-2}}\,K^2}\Big| \leq \frac{\delta}{10K} \text{ for all sufficiently large $i$}. 
\] 
It follows that, for $K$ sufficiently large, 
\[ 
\Big|\partial_r \ln v_i(\hat x_i) + \frac{(n-2)}{r_i}\Big| \leq \frac{\delta}{r_i} \text{ for all sufficiently large $i$}. 
\] 
This contradicts \eqref{Eq:mu=1SSCons}. We conclude that $\bar\beta = \frac{1}{C}$.

\medskip 
\noindent\underline{Step 6:} We wrap up the proof by using the Bishop-Gromov comparison theorem (see e.g. \cite[Theorem 3.10]{ChavelModernIntro}) as in \cite{GV07,TW09}. 
 
We have, by \eqref{RicLB}, 
\[ 
\Ric_{g_{v_i}} \geq -(n-1)\,\epsilon_i^2\,g_{v_i}. 
\] 
where $\epsilon_i = u_i(p)^{\frac{2}{n-2}}\,\alpha \rightarrow 0$. (In the present case, i.e. $\muGp \leq 1$, we can take $\epsilon_i = 0$.)

Let $B^i_r$ and $B^\infty_r$ denote the geodesic balls centered at some point $q \in M \setminus \{x_\infty\}$ (which is fixed for the moment) and of radius $R$ with respect to the metrics $g_{v_i}$ and $g_{v_\infty}$, respectively. By the Bishop-Gromov comparison theorem,  
\[ 
\mu_i(r) := \frac{\textrm{Vol}_{v_i^{\frac{4}{n-2}}g}(B^i_r)}{\mathrm{Vol}(B(\HH^n(-\epsilon_i^2),r)} \text{ is decreasing in r}, 
\] 
where 
\[ 
\mathrm{Vol}(B(\HH^n(-\epsilon_i^2),r) = n\,c(n)\int_0^r \Big(\frac{\sinh(\epsilon_i t)}{\epsilon_i}\Big)^{n - 1}\,dt. 
\] 
Here $c(n)$ is the area of the unit ball in $\RR^n$. Sending $i$ to infinity, we obtain 
\[ 
\mu_\infty(r) := \frac{\mathrm{Vol}_{g_{v_\infty}}(B^\infty_r)}{c(n)\,r^n} \text{ is decreasing in r}. 
\] 

It is clear that $\lim_{r \rightarrow 0} \mu_\infty(r) = 1$. On the other hand, by Step 4 and Step 5, the metric $g_{v_\infty}$ on $M\setminus \{x_\infty\}$ has an `asymptotically flat end at $x_\infty$', in particular, $\lim_{r \rightarrow \infty} \mu_\infty(r) = 1$. It follows that $\mu_\infty \equiv 1$.  
  
Noting that $v_\infty$ is $C^{1,1}_{\rm loc}(M \setminus \{x_\infty\})$, we can proceed as in \cite[Section 3.3]{TW09} to show that $v_\infty$ is smooth.  
 
On the other hand, since $\Ric_{g_{v_i}} \geq -(n-1)\,\epsilon_i^2\,g_{v_i}$ and $v_i$ converges to $v_\infty$ in $C^{1}_{\rm loc}(M \setminus \{x_\infty\})$, we can prove that $\Ric_{g_{v_\infty}} \geq 0$ in $M \setminus \{x_\infty\}$. Indeed, if $\Ric_{g_{v_\infty}}(x_0) < 0$ for some $x_0 \in M \setminus \{x_\infty\}$, then we can find a neighborhood $U$ of $x_0$ in $M \setminus \{x_\infty\}$ and a smooth vector field $X$ supported in $U$ so that $\Ric_{g_{v_\infty}}(X,X) \leq 0$ and $\Ric_{g_{v_\infty}}(X,X)|_{x_0}< 0$ . It follows that 
\begin{align*} 
0  
	&> \int_U \Ric_{g_{v_\infty}}(X,X)\,dv_{g} \\ 
	&= \int_U \Big(-2v_\infty^{-1}\nabla_{g}^2 v_\infty - \frac{2}{n-2} v_\infty^{-1} \Delta_g v_\infty\,g\\ 
		&\qquad\qquad + \frac{2n}{n-2} v_\infty^{-2} dv_\infty \otimes dv_\infty - \frac{2}{n-2} v_\infty^{-2}|dv_\infty|_g^2\,g\Big)(X,X)\,dv_{g}\\ 
	&= \int_U \Big(2v_\infty^{-1}\,[{\rm div}_{g} X\,X(v_\infty) + \nabla_X X(v_\infty)] + \frac{2}{n-2}\,v_\infty^{-1} \nabla_g v_\infty(|X|_g^2)\Big)\,dv_{g}\,\\ 
		&\qquad\qquad + \int_U \Big(\frac{4}{n-2} v_\infty^{-2} dv_\infty \otimes dv_\infty - \frac{4}{n-2} v_\infty^{-2}|dv_\infty|_g^2\,g\Big)(X,X)\,dv_{g}\\ 
	&= \lim_{i \rightarrow \infty} \int_U \Big(2v_i^{-1}\,[{\rm div}_{g} X\,X(v_i) + \nabla_X X(v_i)] + \frac{2}{n-2}\,v_i^{-1} \nabla_g v_i(|X|_g^2)\Big)\,dv_{g}\,\\ 
		&\qquad\qquad + \int_U \Big(\frac{4}{n-2} v_i^{-2} dv_i \otimes dv_i - \frac{4}{n-2} v_i^{-2}|dv_i|_g^2\,g\Big)(X,X)\,dv_{g}\\ 
	&= \lim_{i \rightarrow \infty}\int_U \Ric_{g_{v_i}}(X,X)\,dv_{g}\\ 
	&\geq -(n-1)\epsilon_i^2 \lim_{i \rightarrow \infty}\int_U v_i^{\frac{4}{n-2}}\,g(X,X)\,dv_{g} = 0, 
\end{align*} 
which is absurd. We conclude that $\Ric_{g_{v_\infty}} \geq 0$ on $M$.
 
 We can then invoke the rigidity part of the Bishop-Gromov comparison theorem to obtain that $(M\setminus \{x_\infty\},v_\infty^{\frac{4}{n-2}}g)$ is isometric to $(\RR^n,g_{\rm flat})$. As in \cite[Section 7.6]{GV07}, this implies that $(M,g)$ is conformally equivalent to the standard sphere, contradicting our hypothesis. We have shown that $u \leq C$ in $M$ in the case $\muGp \leq 1$.

 
\subsection{The case $\muGp > 1$} 
 
The argument in Subsection \ref{musmall} carries over except for Step 5. We provide the details for this step in the present case.

In normal coordinates at $x_\infty$, let $r = |x|$. We will need the following lemma, whose proof can be found in Appendix \ref{App:SS}. 
 
\begin{lemma}\label{Lem:mu>1SupSol} 
Assume that $\Gamma$ satisfies \eqref{2}, \eqref{3} and that $\muGp > 1$. For every $1 < \mu < \min(\muGp,2)$ and $0 < \delta < 1$, there is some small $r_1 > 0$ depending only on $(M,g)$, $\mu$ and $\delta$ such that for all $0 < \epsilon < 1$ , the function $v_{\epsilon} = (\epsilon\,r^{-\mu + 1} + 1 - r^\delta)^{\frac{n - 2}{\mu - 1}}$ satisfies 
\[ 
\lambda(A_{g_{v_{\epsilon}}}) \in \Gamma \text{ in } \{0 < r < r_1\}. 
\] 
\end{lemma} 
 
Arguing by contradiction, assume that $a = 0$ in \eqref{vIAsympPrecise}. 
 
Fix $r_1$, $\delta > 0$ and $\mu$ be as in Lemma \ref{Lem:mu>1SupSol}. Decreasing $r_1$ if necessary, we can assume that $r_1 < 1/2$. Let
\[ 
\tilde v_{\epsilon} = 2^{\frac{n-2}{\mu - 1}}\,\Big(\max_{r = r_1} v_\infty\Big)\,v_\epsilon.
\] 
Then $\tilde v_{\epsilon} > v_\infty$ on $\{r = r_1\}$. On the other hand, since $a = 0$ in \eqref{vIAsympPrecise}, there exists $r_2 = r_2(\epsilon) < r_1$ such that $\tilde v_{\epsilon} > v_\infty$ on $\{0 < r \leq r_2\}$. We claim that $\tilde v_{\epsilon} > v_\infty$ on $\{r_2 < r < r_1\}$. If not, then there is some $\theta \geq 1$ such that 
\[ 
\theta\,\tilde v_{\epsilon} \geq v_\infty \text{ in } \{r_2 < r < r_1\} \text{ with equality holds at some $x_0$ with $r_2 < |x| < r_1$}. 
\] 
Since $v_\infty$ satisfies \eqref{9Degenerate} in the viscosity sense and $\theta\,\tilde v_{\epsilon}$ is smooth near $x_0$, this implies that $\lambda(A_{(\theta\,\tilde v_{\epsilon})^{\frac{4}{n-2}}\,g}(x_0)) \in \RR^n \setminus \Gamma$. It is clear from Lemma \ref{Lem:mu>1SupSol} that $\lambda(A_{(\theta\,\tilde v_{\epsilon})^{\frac{4}{n-2}}\,g}(x_0)) \in \Gamma$, a contradiction. The claim is proved and we have $\tilde v_{\epsilon} > v_\infty$ in $\{0 < r < r_1\}$. Sending $\epsilon \rightarrow 0$ we obtain 
\begin{equation} 
v_\infty(x) \leq 2^{\frac{n-2}{\mu - 1}}\,\max_{r = r_1} v_\infty \text{ in } 0 < r < r_1.
	\label{Eq:AlmostRemov} 
\end{equation} 

To complete this step, we need to use a generalization of L\'evy's isopermetric inequality due to B\'erard, Besson and Gallot \cite{BBG85-InvM}. Fix some small $\rho > 0$ for the moment and let $M_{i,\rho} = \{x \in M: d_g(x,x_i) \geq \rho\}$. As seen in Step 1, $\textrm{Vol}_{g_i}(M_{i,\rho}) \rightarrow 0$ as $i \rightarrow \infty$. Thus by \eqref{RicLB}, \eqref{diamBnd} and the isoperimetric inequality \cite{BBG85-InvM}, 
\begin{multline*} 
\left(\int_{\partial B_{g}(x_i,\rho)} u_i^{\frac{2(n-1)}{n-2}}\,dv_g\right)^{\frac{n}{n-1}} 
	= \textrm{Area}_{g_i}(\partial M_{i,\rho})^{\frac{n}{n-1}} \\ 
	\geq C^{-1}\textrm{Vol}_{g_i}(M_{i,\rho})  
	= C^{-1}\int_{M_{i,\rho}} u_i^{\frac{2n}{n-2}}\,dv_g, 
\end{multline*} 
where $C$ is independent of $i$ and $\rho$. Dividing both sides by $u_i(p)^{\frac{2n}{n-2}}$ then sending $i \rightarrow \infty$ we obtain 
\[ 
\left(\int_{\partial B_{g}(x_\infty,\rho)} v_\infty^{\frac{2(n-1)}{n-2}}\,dv_g\right)^{\frac{n}{n-1}} 
	\geq C^{-1}\int_{M \setminus B_{g}(x_\infty,\rho)} v_\infty^{\frac{2n}{n-2}}\,dv_g. 
\] 
For $\rho$ sufficiently small, this cannot happen in view of \eqref{Eq:AlmostRemov} and the positivity of $v_\infty$. This contradiction proves that $a > 0$. 
 
We can now apply the argument in Subsection \ref{musmall} to arrive at a contradiction, which completes the proof of Theorem \ref{Thm:Main}. 
\eproof 
 
Before concluding the section, we note that the argument leading to \eqref{Eq:AlmostRemov} gives the following result, which is of independent interest. 
 
\begin{theorem}\label{Thm:WeakRemSing}
Let $\Omega$ be an open subset of a smooth complete Riemannian manifold $(M,g)$ and $p_0$ be a point in $\Omega$. Assume that $u \in USC(\Omega \setminus \{p_0\} \cap L^\infty_{\rm loc}(\Omega \setminus \{p_0\})$ is a positive viscosity solution  of 
\begin{equation} 
\lambda(A_{g_u}) \in \RR^n \setminus \Gamma \text{ in } \Omega \setminus \{p_0\} 
	\label{9Degenerate} 
\end{equation} 
for some cone $\Gamma$ satisfying \eqref{2}-\eqref{3}. If 
\[ 
\muGp > 1 \text{ and }\limsup_{x \rightarrow p_0} u(x)\,d_g(x,p_0)^{n-2} = 0, 
\] 
then $u$ is locally bounded near $p_0$.  
\end{theorem} 
 
\section{A generalization of Theorem \ref{Thm:Main} when $\muGp \leq 1$} 
\label{general}
 
In this section, we restrict our study to the case $\muGp \leq 1$. In this case, we strengthen Theorem \ref{Thm:Main} to a compactness result for  
\begin{equation} 
f(\lambda(A_{g_u})) =\psi(x)\,u^{-s}, \quad \lambda (A_{g_u}) \in \Gamma,
\quad \mbox{ on } M,
\label{9'new}
\end{equation}
where $s \in [0,\frac{4}{n-2})$ and $\psi$ is a  given positive  smooth function  on $M$.

\begin{theorem}\label{Thm:Cptnessmu<1} 
Let $(f, \Gamma)$ satisfy
(\ref{2})-(\ref{6}) and $\muGp \leq 1$, $(M,g)$ be
a compact, smooth
Riemannian manifold of dimension $n\ge 3$, and $\psi$
be a positive smooth function on $M$. Assume that $(M,g)$ is not conformally equivalent to the standard sphere. For any $\bar s \in [0,\frac{4}{n-2})$, there exists $C = C(M,g,f,\Gamma,\psi,\bar s) > 0$ such that, for any $0 \leq s \leq \bar s$, either \eqref{9'new} has no solution or any positive solution $u \in C^2(M)$ of \eqref{9'new} must satisfy
\[ 
\|\ln u\|_{ C^5(M,g)}
 \leq C.
\] 
\end{theorem}  
 
\bproof The proof follows closely that of Theorem \ref{Thm:Main}. We will only highlight the key changes. Again, for simplicity we only consider $\psi \equiv 1$. Note that local first derivative estimates for \eqref{9'new} are provided by \cite[Theorem 1.10 and Remark 1.12]{Li09-CPAM} while local second derivative estimates for \eqref{9'new} are provided by \cite[Remark 1.13]{Li09-CPAM} and the proof of \cite[Eq. (1.39)]{LiLi03}. Here we have used the assumption $\bar s < \frac{4}{n-2}$.

We will only prove (\ref{maxestimate}) with $\mathfrak{S}_\alpha$ replaced by the solution set of \eqref{9'new}. Arguing by contradiction, assume that there is a sequence $u_i$ of smooth positive functions on $M$, $0 \leq s_i \leq \bar s < \frac{4}{n-2}$ and $x_i \in M$, $x_i \rightarrow x_\infty \in M$ such that the metrics $g_i = u_i^{\frac{4}{n - 2}}\,g$ satisfy 
\[
f(\lambda(A_{g_i})) = u_i^{-s_i}, \qquad \lambda(A_{g_i}) \in \Gamma \;\text{ on } M,
\]
but 
\[
M_i := u_i(x_i) = \max_M u_i \rightarrow \infty. 
\]
 
Note that, since $\muGp \leq 1$, $\Ric_{g_i} \geq 0$.

\medskip 
\noindent \underline{Step 1:} We show that $s_i \rightarrow 0$ and
\begin{equation} 
u_i(x)  \leq Cd_g(x,x_i)^{-\frac{2}{p_i - 1}} \text{ for all } x \in M \setminus \{x_i\},
	\label{IBPXmu<1} 
\end{equation} 
where $p_i = \frac{n+2}{n-2} - s_i$. This is established using the following lemma, which is a variant of Lemma \ref{Lem:NearBlowupPt}.

\begin{lemma}\label{Lem:NBPXmu<1} 
Assume for some $C_1 \geq 1$, $K_i \rightarrow \infty$ and $y_i \in M$ that 
\[
u_i(y_i) \rightarrow \infty \text{ and } \sup \Big\{u_i(y): d_g(y,y_i) \leq K_i\,u_i(y_i)^{-\frac{p_i - 1}{2}}\Big\} \leq C_1\,u_i(y_i) .
\] 
Then $s_i \rightarrow 0$, and for any $0 < \mu < 1$, there exists $K = K(C_1,\mu)$ such that 
\[ 
\textrm{Vol}_{g_i}\Big(\Big\{y: d_g(y,y_i) \leq Ku_i(y_i)^{-\frac{p_i-1}{2}}\Big\}\Big) \geq (1 - \mu) \textrm{Vol}_{g_i}(M). 
\] 
\end{lemma} 
 
\bproof We adapt the proof of Lemma \ref{Lem:NearBlowupPt}. Define $\Phi_i: \RR^n \approx T_{y_i}(M,g) \rightarrow M$ by 
\[ 
\Phi_i(x) = \exp_{y_i} \frac{c^{\frac{p_i - 1}{2}}\,x}{u_i(y_i)^{\frac{p_i - 1}{2}}},
\] 
and let  
\[ 
\tilde u_i(x) = \frac{c}{u_i(y_i)}u_i \circ \Phi_i(x), \qquad x  \in \RR^n. 
\] 
Then $\tilde u_i$ satisfies 
\[ 
f(\lambda(A_{\tilde u_i^{\frac{4}{n-2}}\,\tilde h_i})) = \tilde u_i^{-s_i} \text{ and } \lambda(A_{\tilde u_i^{\frac{4}{n-2}}\,\tilde h_i})) \in \Gamma \text{ on } \{|x| < \delta_0\,c^{-\frac{p_i - 1}{2}}\,u_i(y_i)^{\frac{p_i - 1}{2}}\}, 
\] 
where $\tilde h_i := \Big(\frac{c}{u_i(y_i)}\Big)^{1 - p_i}\Phi_i^* g$ and $\delta_0$ is the injectivity radius of $(M,g)$. As in the proof of Lemma \ref{Lem:NearBlowupPt}, we can assume that $\tilde u_i$ converges in $C^{2,\alpha}_{\rm loc}(\RR^n)$ to some positive $\tilde u_* \in C^2(\RR^n)$ which satisfies 
\[ 
f(\lambda(A_{\tilde u_*^{\frac{4}{n-2}}\,g_{\rm flat}})) = \tilde u_*^{-s_*} \text{ and } \lambda(A_{\tilde u_*^{\frac{4}{n-2}}\,g_{\rm flat}})) \in \Gamma \text{ on } \RR^n,
\] 
where $s_* = \lim_{i \rightarrow \infty} s_i \in [0,\frac{4}{n-2})$. By the Liouville theorem \cite[Theorem 1.3]{LiLi05}, we have $s_* = 0$ and $\tilde u_* = U_{a_*,x_*}$  
for some $a_* > 0$ and $x_* \in \RR^n$ satisfying
\[
|x_*| \leq C \text{ and } C^{-1} \leq a_* \leq C. 
\]
In particular, for any $R > 0$ and $\mu > 0$,  
\begin{equation} 
\|\tilde u_i - \tilde u_*\|_{C^2(\bar B_R)} \leq \mu \text{ for all sufficiently large $i$}. 
	\label{Eq:LiXmu<1} 
\end{equation} 

It follows that the metric $\tilde u_i^{\frac{4}{n-2}}\,\tilde h_i$ converge on compact subsets to the metric $\tilde u_*^{\frac{4}{n-2}}g_{\rm flat}$. Since $(B(0,r), \tilde u_i^{\frac{4}{n-2}}\,\tilde h_i)$ is isometric to $(\Phi_i(B(0,r)),(\frac{c}{u_i(y_i)})^{s_i} g_i)$ for any $r > 0$, we obtain: For any $\epsilon > 0$, there exists $R = R(\epsilon,C_1) > 0$ such that  
\begin{enumerate}[(i)] 
\item $|\textrm{Vol}_{(\frac{c}{u_i(y_i)})^{s_i}g_i}(\Phi_i(B(0,R))) - \textrm{Vol}(\SSphere^n, g_{\rm standard})| \leq C\,\epsilon^n$ for some $C$ independent of $i$ and $\epsilon$, 
\item and the mean curvature of the hypersurface $\partial \Phi_i(B(0,R))$ with respect to $(\frac{c}{u_i(y_i)})^{s_i}g_i$ and the unit normal pointing away from $\Phi_i(B(0,R))$ is no smaller than $\frac{1}{\epsilon}$. 
\end{enumerate} 
Noting that $\Ric_{g_i} \geq 0$, we can apply Proposition \ref{Prop:Hawking} to obtain
\[ 
\textrm{diam}_{(\frac{c}{u_i(y_i)})^{s_i}g_i}(M\setminus\Phi_i(B(0,R))) \leq C\,\epsilon. 
\] 
Thus, by Bishop's theorem, this implies that 
\[ 
\textrm{Vol}_{(\frac{c}{u_i(y_i)})^{s_i}g_i}(M\setminus\Phi_i(B(0,R))) \leq C\,\epsilon^n. 
\] 
Lemma \ref{Lem:NBPXmu<1} is established. 
\eproof   

\medskip 
\noindent\underline{Step 2:} We prove that 
\begin{equation} 
|\nabla_g^k \ln u_i(x)| \leq C\,d_g(x,x_i)^{-k} \text{ for } x \neq x_i, k = 1,2. 
	\label{SC2Estmu<1} 
\end{equation} 
The proof of \eqref{SC2Estmu<1} is exactly as before, except that, instead of \eqref{Eq:wiDef}, we define $w_i$ by 
\[ 
w_i(y) := d_i^{\frac{2}{p_i - 1}}\,u_i(\Psi_i(y)) .
\] 

By Step 1 and Step 2, $u_i$ converges uniformly away from $x_\infty$ to some non-negative limit $u_\infty$, which is either identically zero or always positive. Also, by \eqref{IBPXmu<1}, 
\begin{equation} 
u_\infty(x) \leq C\,d_g(x,x_\infty)^{-\frac{n-2}{2}}. 
	\label{Eq:uinftyBnd} 
\end{equation} 

\medskip 
\noindent\underline{Step 3:} We show that $u_i$ converges in $C^0_{\rm loc}(M\setminus\{x_\infty\})$ to $0$.  
 
The previous argument no longer works. We instead recycle the proof of \eqref{Eq:mu=1viLB}. 
 
Arguing by contradiciton, we assume that the conclusion of Step 3 is incorrect, i.e. $u_\infty > 0$ on $M \setminus \{x_\infty\}$. We will show that there exists some $C > 0$ and $r_1 > 0$ such that for all sufficiently small $\delta > 0$, there holds for some $K = K(\delta) > 0$ and $N = N(\delta) > 0$ that
\begin{equation} 
u_i(x) \geq \frac{1}{C\,d_g(x,x_i)^{n-2-2\delta}} \text{ in } \{K\,u_i(x_i)^{-\frac{p_i - 1}{2}} \leq d_g(x,x_i) \leq r_1\} \text{ for all } i \geq N. 
	\label{Eq:uiLB} 
\end{equation} 

In normal coordinates at $x_i$, let $r = |x|$. Recall the function $\underline{v}_\delta$ and the constant $r_1$ defined in Lemma \ref{Lem:mu=1SubSol}. Since $u_i$ locally converges uniformly away from $x_\infty$ to $u_\infty$ and $u_\infty > 0$ on $M \setminus \{0\}$, there exists some $C$ independent of $i$ such that 
\[ 
u_i \geq \frac{1}{C}\underline{v}_\delta \text{ on } \{r = r_1\}. 
\] 
For some $K = K(\delta) > 0$ to be fixed, let 
\[ 
\bar\beta = \sup\Big\{0 < \beta < \frac{1}{C}: u_i \geq \beta \underline{v}_\delta \text{ in } \{r_i := K\,u_i(x_i)^{-\frac{p_i - 1}{2}} < r < r_1\}\Big\}. 
\] 
We will show that $\bar\beta = \frac{1}{C}$. Assume otherwise that $\bar\beta < \frac{1}{C}$. Then as in the proof of \eqref{Eq:mu=1viLB}, we can find $\hat x_i$ with $|\hat x_i| = r_i$ such that 
\begin{equation} 
\partial_r \ln u_i(\hat x_i) > - \frac{n - 2 - \delta}{r_i}. 
	\label{Eq:uiX1} 
\end{equation} 

On the other hand, by \eqref{Eq:LiXmu<1} with $y_i = x_i$ (note that $\tilde u_* = U_{a,p}$ with $a = 1$ and $p = 0$), 
\[ 
\Big|\frac{c^{\frac{p_i - 1}{2}}}{u_i(x_i)^{\frac{p_i - 1}{2}}}\partial_r \ln u_i(\hat x_i) + \frac{(n-2)c^{-\frac{2}{n-2}}\,K}{1 + c^{-\frac{4}{n-2}}\,K^2}\Big| \leq \frac{\delta}{10K} \text{ for all sufficiently large $i$}. 
\] 
As $p_i \rightarrow \frac{n+2}{n-2}$, it follows that, for $K$ sufficiently large, 
\[ 
\Big|\partial_r \ln u_i(\hat x_i) + \frac{(n-2)}{r_i}\Big| \leq \frac{\delta}{r_i} \text{ for all sufficiently large $i$},
\] 
which contradicts \eqref{Eq:uiX1}. We arrive at $\bar\beta = \frac{1}{C}$, and \eqref{Eq:uiLB}. 
 
Sending $i \rightarrow \infty$ and then $\delta \rightarrow 0$ in \eqref{Eq:uiLB}, we obtain 
\[ 
u_\infty \geq \frac{1}{C\,d_g(x,x_\infty)^{n-2}} \text{ in } \{0 < d_g(x,x_\infty) < r_1\}. 
\] 
But this contradicts \eqref{Eq:uinftyBnd}. We conclude that $u_\infty \equiv 0$.
 
Now define $v_i(x) = u_i(p)^{-1}\,u_i(x)$ and $v_\infty$ to be the limit of $v_i$ as in the proof of Theorem \ref{Thm:Main}.
 
\medskip 
\noindent\underline{Step 4:} We show that $v_\infty$ satisfies \eqref{vIAsympPrecise}. The proof of this statement is exactly as before. 

\medskip 
\noindent\underline{Step 5:} We show that $a > 0$. 
 
The previous argument can be adapted to the current case as in Step 3.
 
\medskip 
\noindent\underline{Step 6:} The conclusion of the proof can be drawn as before. 
\eproof

\section{Proof of Theorem \ref{Thm:Existence}}  
\label{degree}

If $(M,g)$ is conformally equivalent to the standard sphere $\SSphere^n$, the conclusion is clear. We thus assume that $(M,g)$ is not conformally equivalent to $\SSphere^n$.

Fix some $0<\alpha<1$.  For $s\in [0, \frac 2{n-2}]$,
let  
 $$F_s[u] = f(\lambda(A_{g_u})) - u^{-s},
$$
where $u\in C^{4,\alpha}(M)$ satisfies
$u>0$ and $\lambda(A_{g_u}) \in \Gamma$ on $M$.
By Theorem \ref{Thm:Cptnessmu<1}, there exists some 
positive constant $C$ such that every solution of
$F_s[u] = 0$, $0\le s\le \frac 2{n-2}$, satisfies
\begin{equation}
\|\ln u\|_{  C^{4,\alpha}(M) } \le \frac C2
\quad\mbox{and}\quad \mbox{dist}(\lambda(A_{g_u}),
\partial \Gamma)\ge \frac 2 C.
\label{A2-1}
\end{equation}
Let
\[ 
\mathcal{O} = \Big\{u \in C^{4,\alpha}(M):\
u>0,\  \|\ln u\|_{C^{4,\alpha}(M)} \leq C, \lambda(A_{g_u}) \in \Gamma, \dist(\lambda(A_{g_u}), \partial\Gamma) > \frac{1}{C}\Big\} 
\] 
Then
 the degree $\deg(F_s,\mathcal{O},0)$ in the sense of \cite{Li-CPDE} is well-defined and is independent of $s \in [0,\frac{2}{n-2}]$. Thus, to conclude the proof, it suffices to show that $\deg(F_{\frac{2}{n-2}},\mathcal{O},0)
=-1$. 
 
Define, as in \cite{LiLi03} (see page 1424 there),
a homotopy connecting
$(f,\Gamma)$ to $(\sigma_1, \Gamma_1)$:  For $0\le t\le 1$,
$$
\Gamma_t:=\{\lambda\in \Bbb R^n\ |\
t\lambda+(1-t)\sigma_1(\lambda)e\in \Gamma\}, \quad
\mbox{where}\ e=(1,\cdots, 1),
$$
and
$$
f_t(\lambda)=f(t\lambda+(1-t)\sigma_1(\lambda)e).
$$
It was proved in \cite{LiLi03} that $(f_t, \Gamma_t)$ 
also satisfies (\ref{2})-(\ref{6}). 

Consider the problems for $0\le t\le 1$: 
\begin{equation}
f_t(\lambda(A_{g_{u}})) = u^{-\frac{2}{n-2}} \text{ and } \lambda(A_{g_{u}}) \in \Gamma_t \text{ on } M. 
	\label{Eq:SubcritEt}
\end{equation}
This is a family of subcritical equations. The argument leading to \eqref{Eq:LiXmu<1} proves that there is a constant $C > 0$ 
independent of $t$ 
such that all solutions $u$ of \eqref{Eq:SubcritEt} satisfy
\[
u \leq C \text{ on } M.
\]
On the other hand, by evaluating \eqref{Eq:SubcritEt} at a maximum point of $u$, we see that
\[
(\max_M u)^{-\frac{2}{n-2}} \geq (\max_M u)^{-\frac{4}{n-2}}\,f_t(\lambda(A_g)),
\]
which, in view of the assumption $\lambda(A_g) \in \Gamma$ on $M$ and the concavity of $f$, implies that
\[
\max_M u \geq [tf(\lambda(A_g)) + (1 - t)\sigma_1(\lambda(A_g))f(e)]^{\frac{n-2}{2}} \geq \frac{1}{C}.
\]
Hence, by \cite[Theorem 1.10 and Remark 1.12]{Li09-CPAM}, \cite[Eq. (1.39)]{LiLi03} and \cite[Remark 1.13]{Li09-CPAM}, Evans-Krylov's theorem and the Schauder theory, all solutions $u$ of \eqref{Eq:SubcritEt} satisfy (\ref{A2-1}) with $\Gamma$ replaced by $\Gamma_t$. 

Let
$$ 
G_t[u] 
        = f_t(\lambda(A_{g_u})) - u^{-\frac{2}{n-2}}. 
$$
By
 increasing $C$ if necessary, the degree $\deg(G_t,\mathcal{O}_t,0)$ is well-defined and is independent of $t \in [0,1]$ where 
\[ 
\mathcal{O}_t = \Big\{u \in C^{4,\alpha}(M):
\ u>0,  \|\ln u\|_{C^{4,\alpha}(M)} \leq C, \lambda(A_{g_u}) \in \Gamma_t, \dist(\lambda(A_{g_u}), \partial\Gamma_t) > \frac{1}{C}\Big\}. 
\] 
Note that $\deg(G_1,\mathcal{O}_1,0) = \deg(F_{\frac{2}{n-2}},\mathcal{O},0)$.  
 
In the rest of the proof, we 
show that
 $\deg(G_0,\mathcal{O}_0,0)=-1$. Note that $G_0[u] = 0$ is equivalent to 
\[ 
-\Delta_g u + c(n)R_g\,u = u^p \text{ on } M, 
\] 
where $p = \frac{n}{n-2} \in (1,\frac{n+2}{n-2})$
and $c(n)= \frac {n-2}{ 4(n-1)}$. 
It follows that if $u$ is a positive solution of $G_0[u] = 0$ belonging
 to $\mathcal{U} := \{u \in C^{4,\alpha}(M):
\ u>0,  \|\ln u\|_{C^{4,\alpha}(M)} \leq C\}$ then $u \in \mathcal{O}_0$. Hence $\deg(G_0,\mathcal{O}_0,0) = \deg(G_0, \mathcal{U},0)$. 
 
For $0 \leq t \leq 1$ and $p_t = (1 - t) + tp$, let 
\[ 
H_t[u] = -\Delta_g u + [(1 - t) + tc(n)R_g]\,u - \Big[\frac{1 - t}{\textrm{Vol}_g(M)}\int_M u^2\,dv_g + t\Big]u^{p_t}. 
\] 
Note that $H_1[u] = G_0[u]$.

\begin{lemma}\label{Lem:H1Apri} 
For $0 \leq t \leq 1$, positive solutions of $H_t[u] = 0$ satisfy 
\[ 
\|\ln u\|_{C^{4,\alpha}(M)} \leq C(M,g). 
\] 
\end{lemma} 
 
\bproof We only need to consider $t \in (0,1]$, since the case $t = 0$ follows from (a) above. Set 
\[ 
s = \frac{1-t}{\textrm{Vol}_g(M)}\int_M u^{2}\,dv_g + t. 
\] 
Then 
\[ 
(-\Delta_g + [(1- t) + t\,c(n)\,R_g])(s^{\frac{1}{p_t-1}} u) = (s^{\frac{1}{p_t-1}}u)^{p_t} \text{ on } M. 
\] 
Since $u$ is positive and $p$ is subcritical, it is well known that the above equation implies 
\begin{equation} 
\frac{1}{C} \leq s^{\frac{1}{p_t-1}} u \leq C \text{ on } M, 
	\label{Eq:L1su} 
\end{equation} 
where here and below $C$ denotes some constant depending only on $(M,g)$.  
 
From \eqref{Eq:L1su}, we obtain 
\[ 
\max_M u \leq C\min_M u. 
\] 
Thus, from the second inequality in \eqref{Eq:L1su} (at a maximum point of $u$) and the definition of $s$, we have 
\begin{align*} 
C  
	&\geq (\max_ M u)^{p_t - 1}\,((1 - t)\min_M u^{2} + t) \\ 
	&\geq C^{-1}(1 - t)\,(\max_M u)^{p_t+1} + t\,(\max_ M u)^{p_t - 1}, 
\end{align*} 
which implies 
\[ 
u \leq C \text{ on } M. 
\] 
Likewise, from the first inequality in \eqref{Eq:L1su} (at a minimum point of $u$), we obtain 
\[ 
u \geq \frac{1}{C} \text{ on } M. 
\] 

The conclusion follows from standard elliptic estimates applied to the equation $H_t[u] = 0$. 
\eproof 
 
By Lemma \ref{Lem:H1Apri}, the degree $\deg(H_t, \mathcal{U},0)$ is well-defined and independent of $t$. To compute $\deg(H_0,\mathcal{U},0)$, we use the following two facts. 
\begin{enumerate}[(a)] 
  \item If $u = \bar u \equiv 1$ is the unique positive solution of $H_0[u] = 0$. 
  \item If $H_0'[\bar u]\varphi = \mu \varphi$ for some $\mu \leq 0$ and some function $\varphi$ not identically zero, then $\mu = -2
 < 0$ and $\varphi$ is constant. 
\end{enumerate} 
Assuming these facts, it follows from \cite[Propositions 2.3 and 2.4]{Li-CPDE} that $\deg(H_0, \mathcal{U}, 0) = -1$, and so 
\[ 
\deg(F_0,\mathcal{O},0) = -1, 
\] 
as desired.
 
For (a), note that if $u$ is a positive solution of $H_0[u] = 0$, 
then $-1 + \frac{1}{\textrm{Vol}_g(M)}\int_M u^{2}\,dv_g$ is the first eigenvalue of $-\Delta_g$ and $u$ is an associated eigenfunction. It follows that $-1 + \frac{1}{\textrm{Vol}_g(M)}\int_M u^{2}\,dv_g = 0$ and $u \equiv 1$. 
 
For (b), assume that $H_0'[\bar u]\varphi = \mu \varphi$ for some $\mu \leq 0$ and $\varphi \not\equiv 0$. Then 
\begin{equation} 
\mu\varphi = -\Delta_g \varphi - \frac{2}{\textrm{Vol}_g(M)}\int_M \varphi\,dv_g. 
	\label{Eq:(b)varphi} 
\end{equation} 
Integrating over $M$, we obtain 
\[ 
\mu\int_M \varphi dv_g = - 2\int_M \varphi\,dv_g. 
\] 
We claim that $\mu = - 2$. Indeed, if not, the above implies that $\int_M \varphi dv_g = 0$ and so \eqref{Eq:(b)varphi} above implies that $\mu$ is an eigenvalue of $-\Delta_g$ and $\varphi$ is an associate eigenfunction. Since $\mu \leq 0$, it follows that $\mu = 0$ and $\varphi$ does not change sign, which contradicts $\int_M \varphi dv_g = 0$. The claim is proved. Returning to \eqref{Eq:(b)varphi}, we obtain 
\[ 
-\Delta_g \Big[\varphi - \frac{1}{\textrm{Vol}_g(M)}\int_M \varphi\,dv_g\Big] = \mu\Big[\varphi - \frac{1}{\textrm{Vol}_g(M)}\int_M \varphi\,dv_g\Big], 
\] 
Since $\mu < 0$, this leads to 
\[ 
\varphi - \frac{1}{\textrm{Vol}_g(M)}\int_M \varphi\,dv_g \equiv 0, 
\] 
which implies that $\varphi$ is constant. We have proved (b). 
\eproof

\section{A second proof of Theorem \ref{Thm:Main}} 
\label{secondproof}
 
We now provide another proof of Theorem \ref{Thm:Main} that does not use Proposition \ref{Prop:Hawking}. Again, we take for simplicity that $\psi \equiv 1$. The proof for general $\psi$ requires only minor modifications. 

We will only prove \eqref{maxestimate}. The proof of \eqref{minestimate} remains the same once \eqref{maxestimate} is established.
 
Arguing by contradiction, assume that, for some $\alpha \geq 0$, there is a sequence $u_i$ of smooth positive functions on $M$ such that the metrics $g_i = u_i^{\frac{4}{n-2}}g$ satisfy equation \eqref{9i} and the Ricci lower bound \eqref{RicLB} but 
\begin{equation} 
\max_M u_i \rightarrow \infty. 
	\label{BlowupMM} 
\end{equation} 

It is a fact that, for any $R > 0$ and $\epsilon > 0$, there exists a positive constant $C$ depending only on $(M,g)$, $(f,\Gamma)$, $R$ and $\epsilon$ such that, for each sufficiently large $i$, there is a set $\mathcal{S}_i = \{x_i^1, \ldots, x_i^{m_i}\} \subset M$ of finitely many local maximum points of $u_i$ such that
\begin{enumerate}[(i)] 
\item $u_i(x)\,d_g(x,\mathcal{S}_i)^{\frac{n-2}{2}}  \leq C$  for all $x \in M$, 
\item the balls $B_g(x_i^j, R\,u_i(x_i^j)^{-\frac{2}{n-2}})$ are disjoint, 
\item in geodesic normal coordinates (with respect to $g$) at $x_i^j$, 
\[ 
\Big\|u_i(x_i^j)^{-1}u_i\Big(u_i(x_i^j)^{-\frac{2}{n-2}} \cdot\Big) - U_{\frac{1}{2},0}\Big\|_{C^2(B_{2R}(0))} \leq \epsilon, 
\] 
where $U_{\frac{1}{2},0}$ is given by \eqref{UapDef},
\item and $\max_{\mathcal{S}_i} u_i \rightarrow \infty$.
\end{enumerate}
This is a consequence of the Liouville theorem \cite[Theorem 1.3]{LiLi05} and local first and second derivative estimates \cite[Theorem 1.10]{Li09-CPAM} and \cite[Theorem 1.20]{LiLi03}, Evans-Krylov's theorem and the Schauder theory. In the case of the classical Yamabe problem, see \cite{Schoen91}.

If there is no ``bubble accumulation'', i.e.  
\begin{equation} 
\min_{1 \leq j \neq k \leq m_i} d_g(x_i^j,x_i^k) \geq \frac{1}{C} \text{ for some $C$ independent of $i$}, 
	\label{Eq:NoBubbleAcc} 
\end{equation} 
the arguments in Steps 2-6 of the first proof of Theorem \ref{Thm:Main} (with Step 3 being replaced by that in the proof of Theorem \ref{Thm:Cptnessmu<1}) apply and give a contradiction. Indeed, in view of \eqref{Eq:NoBubbleAcc}, $\{m_i\}$ is uniformly bounded. Thus, we can assume without loss of generality that $m_i = m$ is independent of $i$ and that, for each $1 \leq j \leq m$, $x_i^j \rightarrow x_\infty^j \in M$. Step 2 of the first proof shows that
\[
|\nabla_g^k \ln u_i(x)| \leq C\,d_g(x,\mathcal{S}_i)^{-k} \text{ for } x \in M \setminus \mathcal{S}_i, k = 1,2.
\]
By (iv), there is some $j_0$ such that $u_i(x_i^{j_0}) \rightarrow \infty$. In addition, Step 3 of the proof of Theorem \ref{Thm:Cptnessmu<1} shows that if $u_i(x_i^j) \rightarrow \infty$ for some $j$, then $u_i \rightarrow 0$ in $C^0_{\rm loc}(B_g(x_\infty^{j}, r_1) \setminus \{x_\infty^{j}\})$ for some sufficiently small $r_1 > 0$  depending only on $(M,g)$ and the constant $C$ in \eqref{Eq:NoBubbleAcc}. Thus, $u_i \rightarrow 0$ in $C^0_{\rm loc}(M\setminus \mathcal{S}_\infty)$. To apply Steps 4-6 of the first proof, we need to show that $\mathcal{S}_i$ are the blow-up points in the sense that 
\begin{equation} 
\min_{1 \leq j \leq m} u_i(x_i^j) \rightarrow \infty. 
	\label{Eq:SBlowup} 
\end{equation} 
This can be seen as follows: By property (ii), $u_i(x_i^j)$ cannot goes to zero. Thus, if $u_i(x_i^j)$ is bounded for some $j$, property (iii) implies that $u_i$ does not go to zero in a fixed neighborhood of $x_\infty^j$, which contradicts the assertion of Step 3 that $u_i$ goes to zero uniformly away from $\{x_\infty^1, \ldots, x_\infty^m\}$. 
 
In the rest of the proof, we show \eqref{Eq:NoBubbleAcc}. Arguing indirectly, assume that  
\[ 
\ell_i := d_g(x_i^1, x_i^2) = \min_{1 \leq j \neq k \leq m_i} d_g(x_i^j,x_i^k) \rightarrow 0 \text{ as }i \rightarrow \infty. 
\] 

Let $\delta_0$ be the injectivity radius of $(M,g)$. Define $\Xi_i: \RR^n \rightarrow M$ by 
\[ 
\Xi_i(y) = \exp_{x_i^1} (\ell_i\,y). 
\] 
On $\{|y| < \delta_0\ell_i^{-1}\}$, define 
\[ 
\hat u_i(y) = \ell_i^{\frac{n-2}{2}} u_i \circ \Xi_i(y) \text{ and } \hat g_i = \ell_i^{-2}\,\Xi_i^* g. 
\] 
Note that $\hat g_i$ converges on compact subsets of $\RR^n$ to the flat metric $g_{\rm flat}$.  
 
In view of the equation satisfied by $u_i$ and its conformal invariance property. 
\[ 
f(\lambda(A_{\hat u_i^{\frac{4}{n-2}}\,\hat g_i})) = 1 \text{ and } \lambda(A_{\hat u_i^{\frac{4}{n-2}}\,\hat g_i}) \in \Gamma \text{ in } \{|x| < \delta_0\ell_i^{-1}\}. 
\] 

Let $y_i^j = \Xi_i^{-1}(x_i^j)$ and $\hat{\mathcal{S}}_i = \{y_i^1, \ldots, y_i^{m_i}\}$. The $\hat{\mathcal{S}}_i$ satisfies properties (i)-(iii) (but relative to $\hat u_i$). In addition, we also have 
\[ 
 |y_i^j - y_i^k| \geq 1 \text{ for } 1 \leq j \neq k \leq m_i, 
\] 
i.e. an analogue of \eqref{Eq:NoBubbleAcc}. Property (iv) is replaced by: there exists some $R_0 > 0$ such that, 
\begin{equation}
\max_{x \in \hat{\mathcal{S}}_i \cap B_{R_0}} \hat u_i \rightarrow \infty.
	\label{Eq:ivhat}
\end{equation}
Indeed, if this is incorrect, by (i) and (iii) $\hat u_i$ is locally uniformly bounded. Thus, by local first and second derivative estimates \cite[Theorem 1.10]{Li09-CPAM} and \cite[Theorem 1.20]{LiLi03}, Evans-Krylov's theorem and the Schauder theory, $\hat u_i$ converges in $C^2_{\rm loc}$ to some positive limit $\hat u_*$, which by the Liouville theorem \cite[Theorem 1.3]{LiLi05} must have exactly one critical point. On the other hand, each $\hat u_i$ has at least two critical points $y_i^1$ and $y_i^2$, which, in view of (ii) and (iii), converges to two distinct critical points of $\hat u_*$, a contradiction. As before (see \eqref{Eq:SBlowup}), \eqref{Eq:ivhat} implies, for any $r > 0$,
\[
\min_{x \in \hat{\mathcal{S}}_i \cap B_{r}} \hat u_i \rightarrow \infty.
\]

It is clear that the set $\cup \hat{\mathcal{S}}_i$ has isolated accumulation points $\hat{\mathcal{S}}_\infty$. In fact, each of its points are of at least unit distance away from the others. Also, $\hat{\mathcal{S}}_\infty$ has at least two points: it contains $0$ and an accumulation point of $y_i^2$ which has unit modulus. Pick $p \in \RR^n \setminus \hat{\mathcal{S}}_\infty$. We can then follow Steps 2-5 of the first proof of Theorem \ref{Thm:Main} (see the discussion following \eqref{Eq:NoBubbleAcc}) to show that $\hat v_i := \hat u_i(p)^{-1}\,\hat u_i$ converges in $C^{1,\alpha}_{\rm loc}(\RR^n \setminus \hat{\mathcal{S}}_\infty)$ (for all $0 < \alpha < 1$) to some $\hat v_\infty \in C^{1,1}_{\rm loc}(\RR^n \setminus \hat{\mathcal{S}}_\infty)$ which satisfies 
\[ 
\lambda\Big(A_{\hat v_\infty^{\frac{4}{n-2}}\,g_{\rm flat}}\Big) \in \partial \Gamma \text{ in } \RR^n \setminus \hat{\mathcal{S}}_\infty 
\] 
and, for each $y_\infty^j \in \hat{\mathcal{S}}_\infty$, 
\[ 
\lim_{y \rightarrow y_\infty^j} |y - y_\infty^j|^{n - 2}\hat v_\infty(y) \in (0,\infty). 
\] 
The argument in Step 6 of the first proof of Theorem \ref{Thm:Main} then shows $\hat{\mathcal{S}}_\infty$ cannot have more than one points, contradiction our earlier conclusion that it has at least two points.
\eproof

 
\appendix 
\section{Constructions of special sub-solutions and super-solutions} 
 
\label{App:SS} 
 
In this appendix, we give the constructions of sub-solutions and super-solutions of \eqref{9Degenerate} which were needed in the body of the paper.

In the proof we will use the following lemma on the continuity of the eigenvalues of symmetric matrices. 
 
\begin{lemma}\label{Lem:EValCont} 
For an $n\times n$ real symmetric matrix $M$, let $\lambda_1(M), \ldots, \lambda_n(M)$ denote its eigenvalues. There exists a constant $C(n) > 0$ such that for any $\varepsilon > 0$ and any two symmetric matrices $M$ and $\tilde M$ satisfying $|M - \tilde M| < \epsilon$, there holds for some permutation $\sigma = \sigma(M,\tilde M)$ that 
\[ 
\sum_{i = 1}^n |\lambda_i(M) - \lambda_{\sigma(i)}(\tilde M)| < C(n)\epsilon. 
\] 
\end{lemma} 
 
\bproof The result is well known. We present a proof for completeness. Without loss of generality we can assume that $M = (m_{ij}) = \diag(\lambda_1(M), \ldots, \lambda_n(M))$. By Gershgorin's circle theorem, the eigenvalues of $\tilde M = (\tilde m_{ij})$ can be arranged so that 
\[ 
|\lambda_i(\tilde M) - \tilde m_{ii}| \leq \sum_{j \neq i} |\tilde m_{ij}| \text{ for all } i = 1, \ldots, n. 
\] 
Since $M$ is diagonal, this implies that 
\[ 
|\lambda_i(\tilde M) - \lambda_i(M)| \leq \sum_{j = 1}^n |\tilde m_{ij} - m_{ij}| < C\,\epsilon \text{ for all } i = 1, \ldots, n. 
\] 
The assertion follows. 
\eproof

We are now ready to give the proof of Lemmas \ref{Lem:mu>1SupSol}
and  \ref{Lem:mu=1SubSol}.
 
\medskip 
\noindent{\bf Proof of Lemma \ref{Lem:mu>1SupSol}.} Let $\mu$ and $\delta$ be fixed as in the lemma, and $r_1$ will be as small as needed in the proof (though it depends only on $(M,g)$, $\mu$ and $\delta$). Throughout the proof, we will use $C$ to denote some positive constant depending only on $(M,g)$, $\mu$ and $\delta$. The constant $0 < \epsilon < 1$ is arbitrary. 
 
The Schouten tensor for the metric $g_{\epsilon} = v_{\epsilon}^{\frac{4}{n-2}}g$ reads 
\begin{align*} 
A_{g_{\epsilon}}  
	&= - \frac{2}{n-2}\,v_{\epsilon}^{-1}\,\nabla^2_g v_{\epsilon} + \frac{2n}{(n-2)^2}v_{\epsilon}^{-2}\,dv_{\epsilon} \otimes dv_{\epsilon} - \frac{2}{(n-2)^2}\,v_{\epsilon}^{-2}\,|dv_{\epsilon}|_g^2\,g + A_g\\ 
	&= \chi_1\,\mathrm{Id} - \chi_2\,\frac{x}{r} \otimes \frac{x}{r} 
		+ \mathrm{error}. 
\end{align*} 
where $\textrm{Id}$ is the identity matrix and 
\begin{align*} 
\chi_1 
	&= -\frac{2}{n-2}\,v_{\epsilon}^{-1}\,\frac{v_{\epsilon}'}{r} 
		- \frac{2}{(n-2)^2}\,v_{\epsilon}^{-2}\,(v_{\epsilon}')^2\\ 
	&= \frac{2[\epsilon(\mu - 1)\,r^{-\mu + 1} + \delta\,r^\delta][(\mu - 1) - (\mu - 1 + \delta)\,r^\delta]}{(\mu - 1)^2\,r^2(\epsilon\,r^{-\mu + 1} + 1 - r^{\delta})^2},\\ 
	&> \frac{\epsilon\,r^{-\mu + 1} + \frac{\delta}{\mu - 1}\,r^\delta}{C\,r^2(\epsilon\,r^{-\mu + 1} + 1)^2} 	 
		> \frac{1}{C\,r^{2-\delta}(\epsilon\,r^{-\mu + 1} + 1)} \\
	&> \frac{1}{Cr^{3-\mu -\delta}} > 0,\\ 
\chi_2 
	&= \frac{2}{n-2}\,v_{\epsilon}^{-1}\,(v_{\epsilon}'' - \frac{v_{\epsilon}'}{r}) - \frac{2n}{(n-2)^2}v_{\epsilon}^{-2}\,(v_{\epsilon}')^2\\ 
	&= \frac{2}{(\mu - 1)^2\,r^2(\epsilon\,r^{-\mu + 1} + 1 - r^{\delta})^2} 
		\Big[\epsilon(\mu - 1)^2(\mu + 1)\,r^{-\mu + 1}\\ 
		&\qquad\qquad - \epsilon(\mu - 1)((\mu + \delta)^2 - 1)\,r^{-\mu + 1 + \delta}\\ 
		&\qquad\qquad -  \delta(\delta - 2)(\mu - 1)\,r^\delta - 2\delta(\mu + \delta - 1)r^{2\delta}\Big], 
\end{align*} 
and 
\[ 
|\mathrm{error}| \leq C(1 + r\,v_{\epsilon}^{-1}\,|v_{\epsilon}'| + r^2\,v_{\epsilon}^{-2}\,|v_{\epsilon}'|^2) \leq C. 
\] 

Note that 
\[ 
\chi_2 - (\mu + 1)\chi_1 
	= -\frac{2\delta(\mu - 1 + \delta)}{(\mu - 1)\,r^{2-\delta}(\epsilon\,r^{-\mu + 1} + 1 - r^{\delta})} < -\frac{1}{Cr^{3-\mu - \delta}}  < 0. 
\] 
Recalling \eqref{muGpProp} and noting that $\mu < \muGp$, we obtain 
\begin{equation} 
{\rm dist}\Big((1 - \frac{\chi_2}{\chi_1}, 1, \ldots, 1), \RR^n \setminus \Gamma\Big) \geq \frac{1}{C} > 0. 
	\label{WellInside} 
\end{equation} 

Since the eigenvalue of $\chi_1\,\delta_{ij} - \chi_2\,\frac{x}{r} \otimes \frac{x}{r}$ with respect to $\delta_{ij}$ are $\chi_1 - \chi_2, \chi_1, \ldots, \chi_1$, we can apply Lemma \ref{Lem:EValCont} to see that the eigenvalues $\lambda_1, \ldots, \lambda_n$ of $A_{g_{\epsilon}}$ with respect to $g_{\epsilon}$ satisfies 
\begin{align*} 
|\lambda_1 - v_{\epsilon}^{-\frac{4}{n-2}}\,(\chi_1 - \chi_2)| + \sum_{i = 2}^n |\lambda_i - v_{\epsilon}^{-\frac{4}{n-2}}\,\chi_1|  
	\leq C\,v_{\epsilon}^{-\frac{4}{n-2}} \leq C\,r^{3-\mu-\delta}\,v_{\epsilon}^{-\frac{4}{n-2}}\,\chi_1. 
\end{align*} 
It follows that 
\[ 
(\lambda_1, \ldots, \lambda_n) = v_{\epsilon}^{-\frac{4}{n-2}}\,\chi_1 \Big(1 - \frac{\chi_2}{\chi_1} + O(r^{3-\mu-\delta}), 1 + O(r^{3-\mu-\delta}), \ldots, 1 + O(r^{3-\mu-\delta})\Big). 
\] 
Recalling \eqref{WellInside} we arrive at the conclusion for sufficiently small $r_1 > 0$. 
\eproof

\medskip 
\noindent{\bf Proof of Lemma \ref{Lem:mu=1SubSol}.}
Throughout the proof, $C$ denotes some positive constant
depending only on $(M,g)$. 
The constant   $r_1>0$ will be as small as needed in the proof 
(though it depends only on $(M,g)$).
The constant $0<\delta<\frac 14$ is arbitrary.

Let $a=n-2-2\delta$.
 A direct computation shows that the Schouten tensor of the
 metric $\underline{g}_\delta := \underline{v}_\delta^{\frac{4}{n-2}}g$ reads 
\begin{align*} 
	\chi_1\delta_{ij} - \chi_2\,\frac{x}{r} \otimes \frac{x}{r} 
		+ \mathrm{error}. 
\end{align*} 
where 
\begin{align*}
\chi_1
&= \frac 2{  (n-2)^2 }
\frac {  (a-r)(2\delta+r)  }{  r^2}>0,\\
\chi_2 
	&= 
- \frac 2{  (n-2)^2 }
\frac {  -4\delta a+ (n-2-4a)+ 2r^2}  {  r^2}
\end{align*} 
and 
\[ 
|\mathrm{error}| \leq C(1 + r\,\underline{v}_\delta^{-1}\,|\underline{v}_\delta'| + r^2\,\underline{v}_\delta^{-2}\,|\underline{v}_\delta'|^2) \leq C \leq 
\frac {Cr^2}  {2\delta+r}\,\chi_1. 
\] 
Note that 
$$
\chi_2 - 2\chi_1 = 
 \frac 2{  (n-2) r}
\ge \frac {r}{  C(2\delta+r)}
\chi_1.
$$
Thus, by Lemma \ref{Lem:EValCont}, the eigenvalues $\lambda_1, \ldots, \lambda_n$ of $A_{\underline{g}_\delta}$ with respect to $\underline{g}_\delta$ satisfies 
\begin{align*} 
|\lambda_1 -  \underline{v}_\delta^{-\frac{4}{n-2}}\,(\chi_1 - \chi_2)| 
+ \sum_{i = 2}^n |\lambda_i - 
 \underline{v}_\delta^{-\frac{4}{n-2}}\,\chi_1|  
	\leq C\, \underline{v}_\delta^{-\frac{4}{n-2}} \leq 
\frac {Cr^{2}}{  2\delta+r}  \,
 \underline{v}_\delta^{-\frac{4}{n-2}}\,\chi_1. 
\end{align*} 
It follows that 
\begin{align*} 
(\lambda_1, \ldots, \lambda_n)  
	&\leq  \underline{v}_\delta^{-\frac{4}{n-2}}\,\chi_1 \Big(- 1
 - \frac{r}{C(2\delta+r)}, 1 +   \frac{Cr^2}{2\delta+r},
 \ldots, 1 +  \frac{Cr^2}{2\delta+r}\Big)\\ 
	&\leq \Big(1 +  \frac{Cr^2}{2\delta+r}\Big)\, 
\underline{v}_\delta^{-\frac{4}{n-2}}\,\chi_1 
\Big(- 1 -  \frac{r}{C(2\delta+r)}, 1, \ldots, 1\Big), 
\end{align*} 
where we have used the smallness of $r_1$. Since $\muGp \leq 1$,
 $(-1 -  \frac{r}{C(2\delta+r)}, 1,
 \ldots, 1)$ lies outside of $\bar\Gamma$ in view of \eqref{muGpProp}. Thus, $\lambda(A_{\underline{g}_\delta})$ also lies outside of $\bar \Gamma$.
\eproof


\newcommand{\noopsort}[1]{}
\providecommand{\bysame}{\leavevmode\hbox to3em{\hrulefill}\thinspace}
\providecommand{\MR}{\relax\ifhmode\unskip\space\fi MR }
\providecommand{\MRhref}[2]{%
  \href{http://www.ams.org/mathscinet-getitem?mr=#1}{#2}
}
\providecommand{\href}[2]{#2}

\end{document}